\documentclass[11pt]{article}

\usepackage{fullpage,graphicx,authblk,ifthen}
\usepackage{titlesec}
\titlespacing*{\subsection}{0pt}{2.25ex plus 1ex minus .2ex}{1ex plus .2ex}
\titlespacing*{\paragraph}{0pt}{2.25ex plus 1ex minus .2ex}{1em}
\usepackage[authoryear]{natbib}
\usepackage{amsmath,amsfonts,amsthm,amssymb}
\newtheorem{theorem}{Theorem}
\newtheorem{proposition}{Proposition}
\newtheorem{corollary}{Corollary}

\newcounter{def}
\newtheorem{definition}[def]{Definition}

\newtheorem{example}{Example}
   
\newcommand{\myproof}[1]{\textit{Proof.}\hspace{7pt}#1\hfill\qed}
\newcommand{\citt}[1]{\citet{#1}}
\newcommand{\citp}[1]{\citep{#1}}

\usepackage{microtype}  %
\usepackage{booktabs}
\usepackage{multirow}

\usepackage{algorithm}
\usepackage[noend]{algpseudocode}
\usepackage[table,xcdraw]{xcolor}
\usepackage{xspace}
\usepackage{enumitem}
\usepackage{caption,subcaption}
\usepackage{grffile}

\usepackage{hyperref}
\hypersetup{
	colorlinks,
	citecolor=blue,
	filecolor=blue,
	linkcolor=blue,
	urlcolor=black
}
\newcommand{\conv}{\operatorname{conv}}

\newcommand{\inter}{\operatorname{int}}

\renewcommand{\Re}{\mathbb{R}}
\newcommand{\Z}{\mathbb{Z}}
\newcommand{\Q}{\mathbb{Q}}

\newcommand{\F}{\mathcal{F}}

\newcommand{\X}{\mathcal{X}}
\newcommand{\C}{\mathcal{C}}
\newcommand{\N}{\mathcal{N}}

\newcommand{\W}{\mathcal{W}}

\renewcommand{\S}{\mathcal{S}}

\newcommand{\J}{Q}
\newcommand{\RR}{\mathcal{R}}

\renewcommand{\P}{\mathcal{P}}

\newcommand{\R}{\mathcal{V}}

\newcommand{\noprint}[1]{}

\renewcommand{\t}[1]{\texttt{#1}}
\newcommand{\MIBS}{\texttt{MibS}}
\newcommand{\nonl}{\renewcommand{\nl}{\let\nl\oldnl}}
\newcommand{\midd}{\;\middle|\;}
\newcommand{\RNum}[1]{\uppercase\expandafter{\romannumeral #1\relax}}
\newcommand{\NPcomplexity}{\ensuremath{\mathsf{NP}}}

\newcommand{\Cis}[1]{\C_{\rm IS}(#1)}
\newcommand{\Cid}[1]{\C_{\rm ID}(#1)}

\newcommand\norm[1]{\left\lVert#1\right\rVert}
\newcommand\brak[1]{\left\{#1\right\}}

\usepackage{amsmath,amsfonts,commands,xcolor,graphicx,float,amssymb}

\newcommand{\argmin}{\operatorname{argmin}}

\newtheorem{assumption}{Assumption}

\usepackage{pgfplots}
\pgfplotsset{compat=1.16}
\usepackage{enumerate} 
\usepackage{subcaption}
\usepackage{dirtytalk}
\allowdisplaybreaks

\IfFileExists{/home/ted/Ref/BibFiles/LehighMIP.bib}
{

}
{\IfFileExists{/Users/feb223/projects/references/groupreferences/BibFiles/LehighMIP.bib}
{

}
{
    
}}

\graphicspath{{images/}}

\newcommand{\ABSTRACT}{
    We consider the central role of improving directions in solution methods
    for mixed integer bilevel linear optimization problems (MIBLPs). Current
    state-of-the-art methods for solving MIBLPs employ the branch-and-cut
    framework originally developed for solving mixed integer linear
    optimization problems. This approach relies on oracles for two
    kinds of subproblems: those for checking whether a candidate pair of
    leader's and follower's decisions is bilevel feasible, and those required
    for generating valid inequalities. Typically, these two types of oracles
    are managed separately, but
    in this work, we explore their close connection and propose a solution
    framework based on solving a single type of subproblem: determining
    whether there exists a so-called \emph{improving feasible direction} for
    the follower's problem. Solution of this subproblem yields information
    that can be used both to check feasibility \emph{and} to generate strong valid
    inequalities. Building on prior works, we expose the foundational
    role of improving directions in enforcing the follower's optimality
    condition and extend a previously known hierarchy of optimality-based
    relaxations to the mixed-integer setting, showing that the associated
    relaxed feasible regions coincide exactly with the closure
    associated with intersection cuts derived from improving directions. Numerical
    results with an implementation using a modified version of the open source
    solver MibS show that this approach can yield practical improvements.

}

\begin{document}

\title{Improving Directions in Mixed Integer Bilevel Linear Optimization}

\author[1]{ Federico Battista \thanks{\texttt{feb223@lehigh.edu} }}
\author[1]{ Ted K. Ralphs \thanks{\texttt{ted@lehigh.edu}}}
\affil[1]{Department of Industrial and Systems Engineering, Lehigh University, Bethlehem, PA, USA}

\maketitle

\begin{abstract}
    \ABSTRACT
\end{abstract}

\section{Introduction}
\label{sec:intro}

Bilevel optimization problems are a class arising from the recasting of
game-theoretic equilibrium problems, such as that of finding a subgame perfect
Nash equilibrium in the classic Stackelberg game, as mathematical optimization
problems. Stackelberg games are two-player sequential games in which the
players, called the \emph{leader} and the \emph{follower}, make one move each,
with each move consisting of deciding the values of a set of associated
decision variables. The leader chooses values for their decision variables
first and then the follower selects an optimal reaction by solving an
optimization problem called the \emph{follower's (reaction) problem}, with
the leader's solution as an input parameter. In this paper, we focus on the
well-studied special case of a general bilevel optimization problem known as a
\emph{mixed integer bilevel linear optimization problem} (MIBLP) 
in which the decision variables of both the leader and the follower are
constrained by linear inequalities and the values of some variables are
required to be integer.

The goal of solving a bilevel optimization problem is to determine the optimal
decision of the leader. Doing so, however, requires understanding the
functional dependence of the follower's reaction on the leader's decision.
From the leader's standpoint, the casting of the problem as a mathematical
optimization problem can be viewed as introducing constraints involving the
so-called \emph{value function} of the follower's problem into a standard
mathematical optimization problem. The value function encodes optimality
conditions for the follower's problem that are parameterized on the leader's
decision.

The introduction of these additional nonconvex constraints is the main reason
for both the theoretical and practical difficulty of solving problems in this
class. In fact, even the simplest case, in which all constraints are linear
and all variables are continuous, known as the \emph{bilevel linear
optimization problem} (BLP), is
strongly \NPcomplexity-hard~\citp{HanJauSavNew92,BucBilevel23}. Such bilevel
optimization problems are usually tackled by reformulating them as
single-level mathematical optimization problems by incorporating optimality
conditions for the follower's problem derived from, e.g., the KKT conditions,
as constraints; see, e.g.,~\citp{DemZemBilevel20, ZarBorZenProNote19}.

MIBLPs are more difficult to solve than BLPs, both in theory and in practice.
From a theoretical computational complexity standpoint, a decision version of
MIBLP is complete for the class $\Sigma_2^P$~\citp{StoPolynomialtime76,
JerPolynomial85}, the second level of the so-called polynomial hierarchy.
Roughly speaking, this means that even if we had a constant-time oracle for
solving \NPcomplexity-complete subproblems (e.g., MILPs), solving MIBLPs would
still remain as difficult as solving MILPs.

Despite the intractability indicated by the problem's worst-case complexity,
the fact that we can reliably solve MILPs of small- to medium-scale in
practice means that practical algorithms based on oracle computations should
be possible. Indeed, the most successful approach for solving MIBLPs to date
is the branch-and-cut algorithm, which takes a solution approach similar to
that for solving MILPs but relies crucially on the solution of MILP
subproblems both for checking feasibility and for generating the valid
inequalities needed to strengthen the weak initial relaxation.

To date, state-of-the-art branch-and-cut solvers~\citp{FisLjuMonSinNew17,
TahRalDeN20} have treated the MILP subproblems arising when
generating cuts as separate from those arising when checking feasibility of a
solution, despite their theoretical equivalence. We argue that this
separation overlooks an opportunity for significant algorithmic gains. By
adopting a more integrated, ``gray-box'' perspective---where these oracles are
not isolated but instead allowed to share internal information---we can
leverage the strong overlap between separation and feasibility checking. We
show that such a shift reduces redundant oracle computations and paves the way
for more effective solution methods.

\subsection{Definitions and Notation}
\label{sec:miblp}
Before outlining the contribution of this paper, we lay out the formal
definitions and notation, as well as briefly review known results used in
the remainder of the paper. In the literature, a number of equivalent ways of
formulating MIBLPs have been presented. Here, we use the so-called
value function formulation:
\begin{equation}\label{eqn:miblp-vf}\tag{MIBLP}
    \min \left\{cx+d^1y\midd x\in X, y\in \P_1(x) \cap \P_2(x)\cap Y,\; 
    d^2y\leq \phi(b^2 - A^2x)\right\},
\end{equation}
where 
\begin{itemize}
    \setlength\itemsep{0.8em}
    \vspace{-8pt}
    \item $X=\Z^{r_1}_+ \times \Re^{n_1-r_1}_+$ and $Y=\Z^{r_2}_+\times \Re^{n_2-r_2}_+$, 
            respectively, reflect the integrality requirements on the values of the 
            leader's and follower's variables;
    \item $\P_1(x) = \left\{y \in \Re_+^{n_2} \midd G^1 y \geq b^1 - A^1 x\right\}$ 
            is the set of values of the followers variables satisfying the 
            \emph{leader's constraints};
    \item $\P_2(x) = \left\{y\in \Re_+^{n_2}\midd G^2y\geq b^2 - A^2 x \right\}$ 
            is the set of values of the follower's variables satisfying the 
            \emph{follower's constraints},
\end{itemize}
while the input data are $c \in \Q^{n_1}$; $d^1, d^2 \in \Q^{n_2}$;
$A^1\in\Q^{m_1\times n_1}$; $G^1\in\Q^{m_1\times n_2}$; $b^1\in\Q^{m_1}$;
$A^2\in\Q^{m_2\times n_1}$; $G^2\in\Q^{m_2\times n_2}$; and $b^2\in\Q^{m_2}$.
As discussed further below, this formulation assumes the follower behaves in
an \emph{optimistic} fashion (see~\citt{DemFoundations02} for discussion of
other formulations).

The function $\phi$ is the aforementioned \emph{value function} of the follower's
problem, defined as
\begin{equation}\label{eqn:phi} \tag{VF}
    \phi(\beta) = \min \left\{d^2y\midd G^2y\geq \beta, y\in Y \right\} \quad
    \forall \beta \in \Re^{m_2}.
\end{equation}
The role of the value function $\phi$ is to encode the optimality conditions
for the follower's problem. Specifically, for $\hat{x}\in
X, \hat{y} \in \P_1(x) \cap \P_2(x)\cap Y$, if
\begin{equation}\label{eq:opt-constraint}\tag{OPT}
d^2 \hat{y}\leq \phi(b^2 - A^2 \hat{x}), 
\end{equation}
then $\hat{y}$ is contained in the \emph{rational reaction set} w.r.t. to
$\hat{x}$, defined formally as 
\begin{equation}
    \label{eqn:react} \tag{RS}
    \RR(\hat{x}) =  \left\{y \in \S(\hat{x}) \midd d^2 y \leq d^2 \bar{y}, \; \forall \bar{y} \in \S(\hat{x})\right\},
\end{equation}
where
$$\S(\hat{x}) = \left\{ y \in
Y \midd \P_1(\hat{x}) \cap \P_2(\hat{x})\right\}$$ is the set of feasible
solutions for the follower's problem, given the leader's decision. 
In this case, we say that $(\hat{x}, \hat{y})$
is \emph{bilevel feasible}. The set
of all bilevel feasible solutions is the \emph{bilevel feasible region},
denoted by
$$\F = \left\{(x, y) \in \Re^{n_1} \times \Re^{n_2} \midd x \in X,
y \in \RR(x)\right\}.$$ 

When $|\RR(\hat{x})| > 1$, the lower-level problem has multiple
optimal solutions and one can make different assumptions about the procedure
for selecting among alternative optima. The formulation~\eqref{eqn:miblp-vf}
that we employ here specifies the aforementioned \emph{optimistic} assumption in
which the follower selects the response that is most favorable for the
leader's objective function.

The standard relaxations of an MIBLP are either the linear optimization problem (LP) 
relaxation with the polyhedral feasible region 
\begin{equation*}
    \P = \left\{(x,y) \in \Re_+^{n_1} \times \Re_+^{n_2} \midd y \in \P_1(x) \cap \P_2(x)
    \right\},
\end{equation*}  
or the MILP relaxation with feasible region
\begin{equation*}
    \S = \P \cap (X \times Y).
\end{equation*} 
Throughout the paper, we make the following standard assumptions. 
\begin{assumption}
    \label{assume:bounded}
    $\P$ is bounded.   
\end{assumption}

\begin{assumption}
    \label{assume:linking}
    All first-level variables with at least one non-zero coefficient in the second-level problem (so-called linking variables) are integer.  
\end{assumption}

The first assumption guarantees the boundedness of~\eqref{eqn:miblp-vf}, but
is made primarily for the sake of presentation and can be relaxed. The second
assumption ensures that the optimal solution value of~\eqref{eqn:miblp-vf} is
attainable whenever the optimal solution value is
finite~\citp{VicSavJudDiscrete96}. Finally, because the validity of some
inequalities we discuss relies on integrality of the input parameters, we make
the following simplifying assumptions in the remainder of the paper.

\begin{assumption}
$A^2x + G^2y - b^2\in\Z^{m_2}$ for all $(x,y)\in\S$ and $d^2 \in \Z^{n_2}$.
\end{assumption}

\begin{example} 
    For the reminder of this paper, we make use of the well-known bilevel 
    problem from~\citt{MooBarMixed90} as a running example. Figure~\ref{fig:moore_bard} shows
    the bilevel feasible region and optimal solution, along with the standard relaxations 
    $\P$ and $\S$. While a two-dimensional example is illustrative, it 
    cannot capture some of the complexities we address in this work. Hence, we also 
    present a three-dimensional bilevel problem, as shown in Figure~\ref{fig:3d}. 
    
    \begin{figure}[ht!]
        \centering
        \begin{subfigure}[b]{0.47\textwidth}
            \centering
            \includegraphics[width=\linewidth]{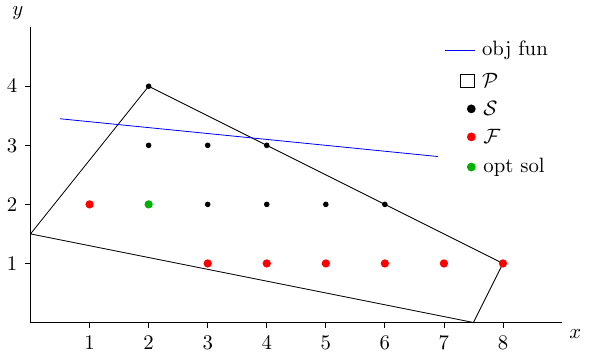} 
        \end{subfigure}
        \hfill
        \begin{subfigure}[b]{0.47\textwidth}
            \centering
            \begin{alignat*}{3}
                \min_{x\in \mathbb{Z}_+} & \quad -x - 10y\notag\\
                \textrm{s.t.} &\quad y\in \operatorname{argmin}\left\{ y:\right.\notag\\
                 & & -5 x + 4 y &\leq 6\notag\\
                 & & x + 2y &\leq 10\\
                 & &  2x - y & \leq 15\notag\\
                 & & 2x + 10y & \geq 15\notag\\
                 & & y &\in \mathbb{Z}_+\left.\right\}\notag
            \end{alignat*}
        \end{subfigure}
        \caption{The feasible region, the LP and MILP relaxation, and the optimal 
                solution of the example from~\citt{MooBarMixed90}.}
        \label{fig:moore_bard}
    \end{figure}
    
    \begin{figure}[ht!]
        \centering
        \begin{subfigure}[b]{0.37\textwidth}
            \centering
            \includegraphics[width=\linewidth]{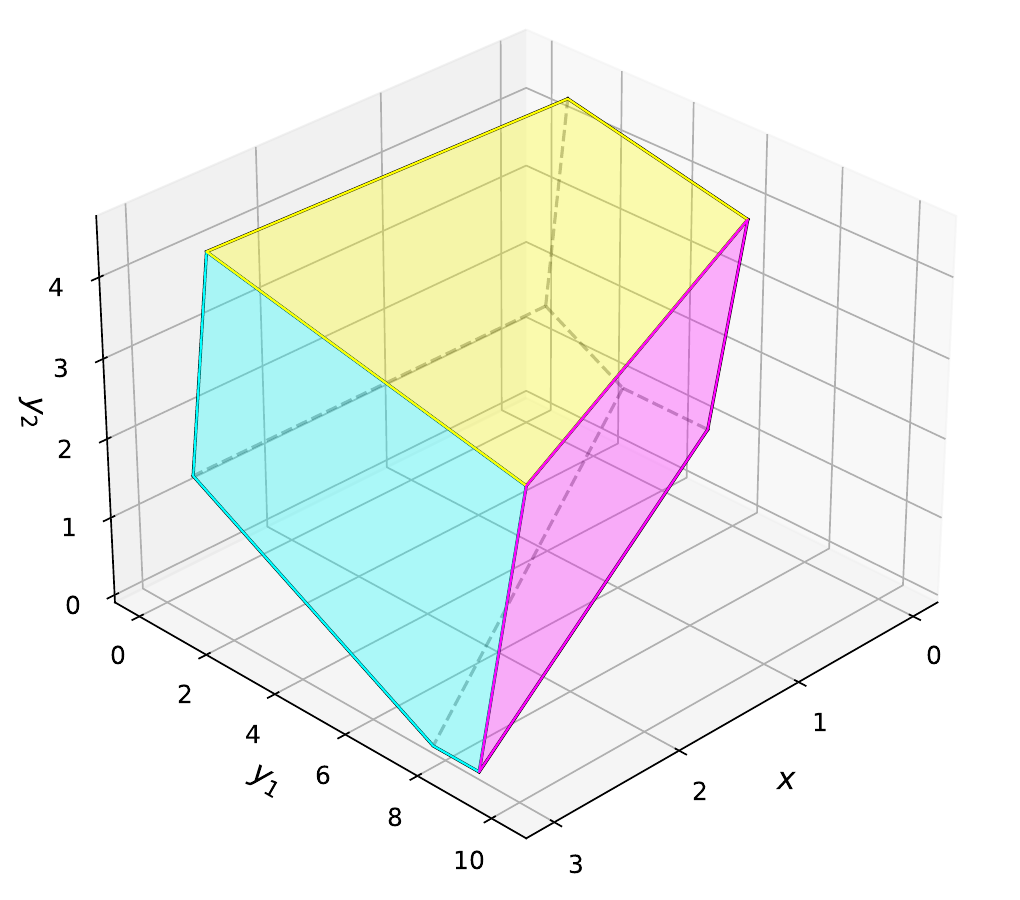} 
        \end{subfigure}
        \hfill
        \begin{subfigure}[b]{0.60\textwidth}
            \centering
            \begin{alignat*}{3}
                \min_{x\in \brak{0, 1, 2, 3}} & \quad -x - 2y_1 - 5y_2\notag\\
                \textrm{s.t.} &\quad y\in \operatorname{argmin}\left\{ y_2:\right.\notag\\
                & & x + 3 y_2 &\geq 3\notag\\
                & & -16x + 12y_1 -4y_2 &\leq 59\label{eqn:3dex}\\
                & &    -x + 9y_1 -2y_2 & \geq 2\notag\\
                & &   - 2x + 6y_1 + 23y_2 & \geq 40\notag\\
                & &   - x + y_1 + 10y_2 & \leq 45\notag\\
                & &   y &\in \mathbb{Z}^2_+\left.\right\}\notag
            \end{alignat*}
        \end{subfigure}
        \caption{A three-dimensional MIBLP and the feasible region of its LP relaxation.}
        \label{fig:3d}
    \end{figure}
\end{example}

\subsection{Improving Directions}
\label{sec:id}

We now introduce the concept of improving directions. 
\begin{definition}
    \label{def:id}
    A vector $w \in \Z^{r_2} \times \Re^{n_2-r_2}$ is an \emph{improving
    direction} (ID) if $d^2w < 0$ and we denote the set of all IDs as $\W
    := \{w \in \Z^{r_2} \times \Re^{n_2-r_2} \mid d^2w < 0\}$.
    With respect to a given $(\hat{x}, \hat{y}) \in \P$, an ID $w$ is an
    \emph{improving feasible direction} (IFD) if $\hat{y} + w \in \P_2(\hat{x})$. The set of all 
    improving feasible directions with respect to $(\hat{x}, \hat{y})$ is 
\begin{equation}
    \label{eqn:fid}
    \W(\hat{x}, \hat{y}) = \left \{ w \in \W \mid 
    \hat{y} + w \in \P_2(\hat{x}) \right \}.
\end{equation}
\end{definition}
Informally, improving feasible directions with respect to a bilevel infeasible
solution are those that point towards the bilevel feasible region $\F$. Note
that although this definition requires that improving directions must
themselves satisfy the integrality requirements of the follower's problem, an
improving direction can nevertheless be said to be feasible with respect to
any point in $\P$, not only points in $\S$.

The most obvious use of improving feasible directions is as a certificate of
bilevel infeasibility for points in $\S$.
We present the following fundamental result that shows the relationship between
the existence of an improving direction and bilevel feasibility of a given point.
\begin{proposition}
    \label{lemma:bfeas}
    Let $(\hat{x}, \hat{y}) \in \S$. Then we have $(\hat{x}, \hat{y}) \in \F \Longleftrightarrow \W(\hat{x}, \hat{y}) = \emptyset$.
\end{proposition}
\myproof{
    ($\Rightarrow$) Assume $(\hat{x}, \hat{y}) \in \F$ and $\W(\hat{x}, \hat{y}) 
    \not = \emptyset$ for sake of contradiction. Then $d^2\hat{y} = \phi(b^2 - A^2\hat{x})$.
    Now let $w \in \W(\hat{x}, \hat{y})$, then by definition $\hat{y} + w \in \P_2(\hat{x})$.
    Moreover, we have
    $$ d^2(\hat{y} + w) < d^2\hat{y} = \phi(b^2 - A^2\hat{x}). $$
    This implies that either $d^2\hat{y} > \phi(b^2 - A^2\hat{x})$ or 
    $\hat{y} \notin \P_2(\hat{x})$, then $\hat{y} \notin \RR(\hat{x})$. 
    This contradicts the bilevel feasibility of $(\hat{x}, \hat{y})$. \\\\
    \noindent
    ($\Leftarrow$) We prove the contrapositive. Let $(\hat{x}, \hat{y}) \notin \F$ be given. Then
    $\exists\ \bar{y} \in \S(\hat{x})$ such that $d^2\bar{y} < d^2\hat{y}$. Now consider 
    $w := \bar{y} - \hat{y}$. Note that $\hat{y} + w \in \P_2(\hat{x})$ and $d^2w < 0$ by construction.
    Then $w \in \W(\hat{x}, \hat{y})$ and the statement is proven.
}   

This result leads to an alternative method to check bilevel feasibility, which
is to check whether $\W(\hat{x}, \hat{y}) = \emptyset$. By
Proposition~\ref{lemma:bfeas}, when $(\hat{x}, \hat{y}) \in \S$, emptiness of
$\W(\hat{x}, \hat{y})$ is equivalent to bilevel feasibility. On the other
hand, $\W(\hat{x}, \hat{y})$ may be empty if
$(\hat{x}, \hat{y}) \in \P \setminus \S$, even if
$(\hat{x}, \hat{y}) \not\in \conv(\F)$ (Proposition~\ref{lemma:bfeas} does not
apply). From a computational perspective, this is quite important. As a side
note, Proposition~\ref{lemma:bfeas} also indicates that checking feasibility of a
given solution is a problem in co-NP, which is interesting, though not
unexpected. 

Given a candidate pair $(\hat{x}, \hat{y}) \in \S$ and a direction 
$w \in \W(\hat{x}, \hat{y})$, the follower's solution can be 
augmented to form an \emph{improving solution}, i.e., 
a new candidate pair $(\hat{x}, \hat{y} + w)$ that remains 
in $\S$ but improves the follower's objective, i.e., 
$d^2(\hat{y} + w) < d^2\hat{y}$.
In fact, when $(\hat{x}, \hat{y})$ satisfies integrality requirements, checking whether 
$\W(\hat{x}, \hat{y})$ is empty is formally equivalent to 
determining the existence of an improving solution to the 
follower's problem, as both can be formulated as MILPs.
However, as discussed in Section~\ref{sec:genid}, elements of 
$\W(\hat{x}, \hat{y})$ can be generated using a variety of objective 
functions, enabling the computation of directions with different 
desirable properties. This makes it a seemingly more flexible 
approach to certifying bilevel infeasibility.

In the context of traditional, single-level integer linear optimization, the
problem of determining whether there exists a so-called \emph{augmenting
vector} has been previously studied and it was shown to be
oracle-polynomial-time equivalent to solving an MILP by~\citt{SchRelative09},
who called it the \emph{augmentation problem}. Similarly, bilevel feasibility
of a given candidate solution in the MIBLP context can be checked by
determining whether an IFD exists. It follows that an algorithm for solving
MIBLPs based only on generating IFDs is possible.

\begin{example}
In Figure~\ref{fig:moore_bard_id} we show three bilevel infeasible points and possible improving feasible 
directions moving them toward points in $\F$ of the~\citt{MooBarMixed90} example. 
Note that for points in $\S \setminus \F$, there must exist
at least one IFD (see, e.g., $\hat{y}^1$ and $\hat{y}^2$). However, IFDs may also
exist for points in $\P$ (e.g., $\hat{y}^3$). Moreover, the same improving direction (e.g., $w_1$) may be 
feasible for more than one point. It is easy to verify that given any of these points (say, $\hat{y}^1$), 
we can move along any of the IFDs (say, $w_1$), to obtain an improving feasible solution.
\begin{figure}[h]
    \centering
    \includegraphics[width=0.6\linewidth]{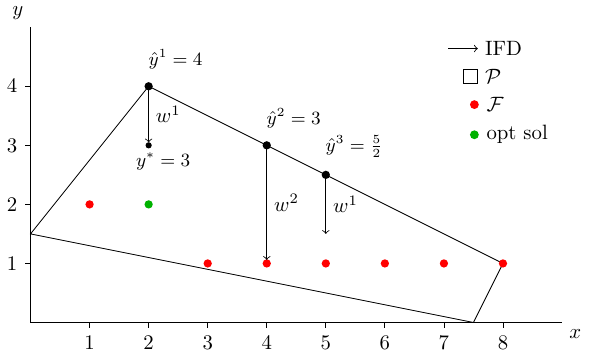}
    \caption{The improving feasible directions moving points in $\S$ (or $\P$) towards points in $\F$.}
    \label{fig:moore_bard_id}
\end{figure}
\end{example}

\subsection{Valid Inequalities}

We briefly review the basic concepts of valid inequalities.

\begin{definition}[Valid Inequality]
A \emph{valid inequality} for $\F$ is a triple
$(\alpha^x, \alpha^y, \beta) \in \Q^{n_1} \times \Q^{n_2} \times \Q$ such that
    $$ \F \subseteq \brak{(x,
    y) \in \Re^{n_1} \times \Re^{n_2} \midd \alpha^xx + \alpha^y 
    y \geq \beta}. $$ 
\end{definition}

The goal of generating valid inequalities in a branch-and-cut algorithm can
be expressed in several related ways. Most obviously, the goal is to
strengthen the relaxation, leading to an improved dual bound. This is
accomplished by removing solutions to the relaxation that are infeasible to
the original problem through the addition of valid inequalities violated by
those solutions.
This shrinks the feasible
region of the relaxation, resulting in a better approximation of $\conv(\F)$.

Another way of viewing the same goal is as that of implicitly enforcing
constraints that were relaxed, in this case the optimality conditions.
Although the constraints that were relaxed are nonlinear, their effect can be
replicated using linear inequalities because the objective function is linear,
which allows us to convexify the feasible region. Thus, while the goal in the
MILP setting is to expose points satisfying integrality conditions as extreme
points of the relaxed feasible region, here we need to additionally ensure
that the exposed solutions satisfy second-level optimality conditions. In
other words, the exposed integer points must also lie on the boundary of the
epigraph of the value function (roughly speaking). The fact that the exposed
points must simultaneously satisfy these two different properties results in
complex interactions when generating valid inequalities that require careful
algorithmic control. More details of the theory underlying the generation of
valid inequalities in the MIBLP case, as well as detailed discussions of
control mechanisms are provided in~\citt{TahRal25}.

The general recipe by which most of the known classes of valid inequalities
for MIBLPs are constructed is to first identify a set $\C$ containing no
(improving) bilevel feasible solutions in its interior and then generate an
inequality valid for $\conv(\overline{\inter(\C)} \cap \P)$. Such an
inequality is valid for $\F$, since
$\F \subseteq \overline{\inter(\C)} \cap \P$. Because this recipe is strictly
a generalization of the one first proposed by~\citt{BalInteger72} in the
context of MILPs, these inequalities are sometimes referred to broadly
as \emph{intersection cuts} (ICs). By employing different ``solution free''
sets and by replacing $\F$ with the feasible region of a relaxation, we can
derive a wide range of different classes of valid inequality.

The most common way in which ICs are generated in practice is the method
proposed in the original paper of~\citt{BalInteger72}. That is, we replace
$\P$ with a simplicial radial cone that contains $\F$ and whose single extreme
point lies in the interior of a convex set $\C$ containing no feasible points in
its interior, as described above. Then the hyperplane defined by the points of
intersection of the rays of the simplicial cone with the set $\C$ is a
hyperplane that separates the cone's extreme point from $\F$. When we
reference the term \emph{intersection cut} in the remainder of the paper, we
are referring to this specific type of IC, defined as follows.
\begin{definition}[Bilevel Free Set]
    A \emph{bilevel free set} (BFS) is a set $\C \subseteq \Re^{n_1 + n_2}$ 
    such that $\inter(\C) \cap \F = \emptyset$. 
\end{definition}

\begin{definition}[Intersection Cut] \label{def:IC}
Let $\C \subseteq \Re^{n_1+n_2}$ be a convex BFS with a given point
$(\hat{x}, \hat{y}) \in \Re^{n_1+n_2}$ in its interior. Let
$\R(\hat{x}, \hat{y})$ be a simplicial radial cone containing $\F$ with vertex
$(\hat{x}, \hat{y})$. Then, if the triple
$(\alpha^x,\alpha^y, \beta)\in\Q^{n_1+n_2+1}$ is such that the set $\{(x,
y) \in \Re^{n_1 + n_2} \mid \alpha^x x + \alpha^y y = \beta\}$ is the unique
hyperplane containing the points of intersection of $\C$ with the extreme rays
of $\R(\hat{x}, \hat{y})$, we have $(\alpha^x,\alpha^y, \beta)$ is an
inequality valid for $\F$ and violated by $(\hat{x}, \hat{y})$. Such
inequality is called an \emph{intersection cut}.
\end{definition}

The simplicial cone is often taken to be that described by a linearly
independent set of inequalities binding at some basic feasible solution to the
initial LP relaxation, but we consider cones derived in other ways later in
the paper. In Section~\ref{sec:IC-closure}, we define a notion of rank
similar to the standard notion from the theory of valid inequalities for
MILPs. Beginning with inequalities derived from (bilevel infeasible) extreme
points of $\P$, which are the inequalities of rank 1, one can then iteratively
apply this procedure to derive inequalities of higher rank.

In~\citt{FisLjuMonSinNew17,FisLjuMonSinUse18}, a number of classes of ICs were
introduced, but two in particular play a central role in what follows---those
arising from the existence of \emph{improving solutions} and \emph{improving
directions}, respectively. In~\citp{TahRal25}, the authors refer to these
two classes as \emph{improving solution intersection cuts} (ISICs) and
\emph{improving direction intersection cuts} (IDICs). 

As observed by~\citt{FisLjuMonSinNew17}, the strength of an IC is directly
related to the \say{size} of the set $\C$: a \say{larger} such set should
result in a stronger inequality (see discussion
in~\cite{FisLjuMonSinNew17,FisLjuMonSinUse18}). Because of the specific forms
of convex sets utilized for the two classes of ICs, an improving
direction/solution that will result in a large bilevel free set for one class
will not necessarily result in a large bilevel free set for the other class.

\subsection{Contribution}

Improving directions arise naturally in bilevel optimization and the use of an
oracle for finding such directions in an algorithm for solving bilevel
optimization problems is not new. As observed in~\citt{TahRal25}, ICs
generated from improving directions are the strongest inequalities known from
an empirical standpoint. Beyond separation, such oracles have also been
used for branching~\citp{WanXuWatermelon17} and for strengthening
relaxations~\citp{XueProRal22}.

The contributions of this paper are two-fold. From the theoretical standpoint,
we highlight the fundamental role of improving directions in restoring the
follower's optimality condition. To this end, we generalize the hierarchy of
optimality-based relaxations presented in~\citt{XueProRal22} for bilevel
linear problems with only binary variables to the mixed-integer context and
show that the convex hulls of the feasible regions defined by this hierarchy
can be exactly characterized using inequalities generated from improving
directions.

From a computational standpoint, our contribution is to show that unifying the
oracle computations for the separation of bilevel infeasible points with those
for checking bilevel feasibility and leveraging the equivalence of the
underlying MILP subproblems can lead to significant improvements in the
practical performance of the branch-and-cut algorithm. Given that the problem
of finding an improving direction serves both to generate strong valid
inequalities and check bilevel feasibility, it naturally arises as a suitable
candidate for the purpose. Motivated by this, we propose and implement a
branch-and-cut algorithm that follows the basic outline described
in~\citt{TahRalDeN20} but avoids checking feasibility of solutions by solving
the follower's problem and instead solves the problem of determining whether
there exists an improving feasible direction. We argue that solving the
leader's problem to provable optimality without ever explicitly evaluating the
follower's value function is not just possible but has the potential to
improve empirical performance.

\subsection{Outline}

The remainder of this paper is organized as follows: 
Section~\ref{sec:k-opt} introduces a hierarchy of optimality-based relaxations
that arise from improving directions. 
Section~\ref{sec:valid} discusses the generation of 
inequalities valid for bilevel feasible points derived from improving directions.
Section~\ref{sec:algo} introduces the algorithm and 
presents methods to generate these improving feasible directions.
Finally, Section~\ref{sec:res} presents the computational results.

\section{The $k$-opt Hierarchy}
\label{sec:k-opt}

Both of the standard relaxations of MIBLPs discard the optimality conditions
for the follower's problem completely, which results in weak relaxations
in general. In zero-sum problems, for example, the follower's problem is
implicitly being solved with an objective function that is precisely the
opposite of the follower's true objective.

Rather than completely discard the optimality condition, it is possible to
impose a weaker condition that improves the bound yielded by the linear
relaxation, yet whose computation remains tractable. In general, this can be
achieved by replacing the reaction set $\RR(x)$ with a suitable superset.
A straightforward way to define such a set is to replace the follower's
problem with a relaxation. 
In~\citt{XieRalPro25}, for example, the optimality condition is
relaxed by allowing lower-level solutions within a specified optimality gap.
Here, we follow the ideas of~\citt{XueProRal22} and instead require the
follower's reactions to satisfy local rather than global optimality
conditions.
In~\citp{XueProRal22}, this was done for MIBLPs with only binary
variables in the follower's problem. We extend the method by
generalizing the definition of \emph{$k$-neighborhood}. For a given
$y \in Y$ and $k \in \Z_+$, the $k$-neighborhood of $y$ is 
$$ \N_k(y) = \left\{ \bar{y} \in Y \midd \norm{\bar{y} - y}_1 \leq k \right\}.$$
For $y \in Y$, the elements in the $k$-neighborhood are the points in $Y$ that
can be reached by following any direction $w \in Y$ with $\norm{w}_1 \leq k$.
Therefore, we can modify the reaction set~\eqref{eqn:react} to one requiring
only local optimality by replacing $\RR(x)$ with
\begin{equation}
    \label{eqn:kopt-react} \tag{$k$-RS}
    \RR(x; k) = \left\{y \in \S(x) \midd d^2 y \leq d^2 \bar{y}, \; \forall \bar{y} \in \N_k(y) \cap \S(x)\right\}.
\end{equation}
for some $k \in \Z_+$. Note that for $k = 0$ we have that $\N_0(y) = \brak{y}$ and $\RR(x; 0) = \S(x)$, whereas for 
$$ \bar{k} := \sum_{i = 1}^{r_2} \left( \left \lfloor \max_{(x, y) \in \P} y_i \right \rfloor  - \left \lceil \min_{(x, y) \in \P} y_i \right \rceil \right) 
+ \sum_{i = r_2 + 1}^{n_2} \left\lceil \max_{(x, y) \in \P} y_i - \min_{(x, y) \in \P} y_i \right\rceil,$$
which is finite by Assumption~\ref{assume:bounded}, we have $\N_{\bar{k}}(y) =
Y$ and $\RR(x; \bar{k}) = \RR(x)$, for all $x \in X$. 

In~\citt{XueProRal22}, \eqref{eqn:kopt-react} is referred to as
the \emph{$k$-optimal reaction set} and any $y \in \RR(x; k)$ as
a \emph{$k$-optimal reaction}. Moreover, the authors show that for $Y
= \left\{ 0, 1 \right\}^{n_2}$ and for any fixed $k \in \Z_+$, the set of
points in \eqref{eqn:kopt-react} is MILP-representable with a description of
polynomial size. 

The MIBLP with $k$-optimal follower is then formally defined as follows:
\begin{align}
    \label{eqn:k-miblp}
           \min &\; cx + d^1y       \notag\\
    \mbox{s.t.} &\; (x, y) \in \S   \tag{$\mbox{BP}_k$}\\
                &\; y \in \RR(x; k), \notag 
\end{align}
for all $k \in \Z_+$. The set of feasible points of~\eqref{eqn:k-miblp} is 
\begin{equation}
    \label{eqn:k-feas}
    \F(k) = \brak{(x, y) \in \Re^{n_1} \times \Re^{n_2} \midd x \in X, y \in \RR(x; k)}.
\end{equation}
Following an approach similar to~\citt{XueProRal22}, we can characterize
membership in $\RR(x; k)$ for $x \in X$ by a set of necessary and sufficient
conditions. For this purpose, we denote the set of all improving
directions with a 1-norm no bigger than $k$ as 
$$\W^k = \left \{ w \in \W \mid \norm{w}_1 \leq k \right \},$$
and those that are feasible with respect to a given $(\hat{x}, \hat{y})$ 
as

$$\W(\hat{x}, \hat{y}; k) = \left \{ w \in \W^k \mid 
    \hat{y} + w \in \P_2(\hat{x}),\ \norm{w}_1 \leq k \right \}.$$
Then we have the following result.
\begin{proposition}
    \label{lemma:k-bfeas}
    Let $(\hat{x}, \hat{y}) \in \S$. Then $(\hat{x}, \hat{y}) \in \F(k) \Longleftrightarrow \W(\hat{x}, \hat{y}; k) = \emptyset$.
\end{proposition}
\myproof{
    Let $k \in [\bar{k}]$. Then the proof is divided into two parts.\\
    ($\Rightarrow$) We prove the contrapositive. Let $\W(\hat{x}, \hat{y}; k) \not = \emptyset$ be given. 
    Then there exists $w \in \W(\hat{x}, \hat{y}; k)$ such that $\bar{y} := \hat{y} + w \in \N_k(\hat{y}) \cap \S(\hat{x})$ 
    with $d^2\bar{y} < d^2\hat{y} \Rightarrow \hat{y} \notin \RR(\hat{x}; k) \Rightarrow (\hat{x}, \hat{y}) \notin \F(k)$. \\\\
    \noindent
    ($\Leftarrow$) Again, we prove the contrapositive. Let $(\hat{x}, \hat{y}) \notin \F(k)$ be given. Then 
    $\hat{y} \notin \RR(\hat{x}; k) \Rightarrow \exists\ \bar{y} \in \N_k(\hat{y}) \cap \S(\hat{x})$ with $d^2\bar{y} < d^2\hat{y}$.
    Now consider $w := \bar{y} - \hat{y}$. Note that $\norm{w}_1 \leq k$, $w \in \P_2(\hat{x})$ and $d^2w < 0$ by construction. Then
    $w \in \W(\hat{x}, \hat{y}; k)$ and this proves the statement.
}

In other words, just as we can check the feasibility of $(\hat{x}, \hat{y})$
with respect to the constraints of~\eqref{eqn:miblp-vf} by checking emptiness
of $\W(\hat{x}, \hat{y})$, we can check feasibility with respect
to~\eqref{eqn:k-miblp} by checking emptiness of $\W(\hat{x}, \hat{y}; k)$.

With the next result, we show that~\eqref{eqn:k-miblp} defines a hierarchy of
relaxations for~\eqref{eqn:miblp-vf}.

\begin{theorem}\label{thm:hierarchy}
    $\S = \F(0) \supseteq \F(1) \supseteq \F(2) \supseteq \ldots \supseteq \F(\bar{k}) = \F.$
\end{theorem}
\myproof{
    To show that $\F(k) \supseteq \F(k + 1)$ for $k = 0, 1, \ldots, \bar{k} - 1$, it is sufficient to observe
    that $\W(x, y; k^\prime) \subseteq \W(x, y; k^{\prime\prime})$, for all $k^\prime, k^{\prime\prime} \in [\bar{k}]$, 
    with $k^\prime \leq k^{\prime\prime}$, then it follows from Proposition~\ref{lemma:k-bfeas}. For $k = 0$ we have 
    that $\N_0(y) = \brak{y}$ and $\RR(x; 0) = \S(x)$. Furthermore, if $k = \bar{k}$ then $\N_{\bar{k}}(y) = Y$ and $\RR(x; \bar{k}) = \RR(x)$.
}

The next question we address explores the theoretical computational complexity of 
computing the dual bound of~\eqref{eqn:k-miblp} for a fixed value of $k$. 
As a natural extension of the result in~\citt{XueProRal22}, the following 
theorem shows that the decision version of~\eqref{eqn:k-miblp} is \NPcomplexity-complete.

\begin{theorem}
    \label{thm:k-opt-complexity}
    The decision version of~\eqref{eqn:k-miblp} is \NPcomplexity-complete for any fixed integer $k \geq 1$.
\end{theorem}
\myproof{
    Given any fixed integer $k \geq 1$, we consider the feasibility problem 
    associated with~\eqref{eqn:k-miblp}, which is one form of decision version
    of ~\eqref{eqn:k-miblp}. That is, we consider the problem of deciding
    whether $\exists (x, y) \in \F(k)$. First, we prove that this decision
    problem is in \NPcomplexity. Then, we show that the problem of determining
    whether a given MILP is feasible can be reduced to that of deciding
    whether~\eqref{eqn:k-miblp} is feasible.

    To show that the feasibility problem associated with~\eqref{eqn:k-miblp}
    is in \NPcomplexity, we show that when $\F(k) \not= \emptyset$, then there
    exists a certificate that can be verified in polynomial time. When
    $\F(k) \not= \emptyset$, there must exist $(\hat{x}, \hat{y}) \in \F(k)$. By Proposition~\ref{lemma:k-bfeas}, $(\hat{x}, \hat{y}) \in \F(k)$
    if and only if $\W(\hat{x}, \hat{y}; k) = \emptyset$. Then, a verifier is
    deciding whether $\forall w \in \W^k, w \notin \W(\hat{x}, \hat{y}; k)$.
    Given any $w \in \W^k$, verifying that $w \in \W(\hat{x}, \hat{y}; k)$ can
    be done in polynomial time by deciding the membership of $\hat{y} + w$ to
    the polytope $\P_2(\hat{x})$. Moreover, the set $\W^k$ has a cardinality
    of $\mathcal{O}(n^k)$, which is polynomial in $n$ for any fixed integer
    $k \geq 1$. Therefore, the certificate can be verified in time polynomial in $n$.

    To show that the problem of deciding feasibility of an MILP can be reduced
    to that of deciding feasibility of~\eqref{eqn:k-miblp}, let $\X
    := \left\{ x \in X \midd A^1x \geq b^1 \right\}$ and consider the MILP
    feasibility problem of determining whether $\X \not= \emptyset$. The
    certificate for this problem is any $x \in \X$. We show how to construct
    such a certificate from the certificate of an instance
    of~\eqref{eqn:k-miblp}. As such, let an instance of~\eqref{eqn:k-miblp}
    feasibility be defined as follows. We let $Y = \Re_+^{n_2}$, $G^1 =
    0_{m_1 \times n_2}$, $A^2 = 0_{m_2 \times n_1}$, $G^2 = I_{n_2}$ (the
    identity matrix of order $n_2$), $b^2 = \mathbf{0}_{n_2}$ and $d^2
    = \mathbf{1}_{n_2}$ (the all-zeros and all-ones vectors of order $n_2$, respectively). 
    Note that by construction, $\RR(x; k) = \RR(x)
    = \{\mathbf{0}_{n_2}\}$ for all $x \in \X$. It is easy to verify that
    $x \in \X \Longleftrightarrow (x, \mathbf{0}_{n_2}) \in \F(k)$. Then any
    certificate of feasibility for the constructed instance of~\eqref{eqn:k-miblp} 
    can be mapped to a certificate for the feasibility of the MILP. 
}

In the context of bilevel problems with only binary variables at the second level, 
\citt{XueProRal22} show that optimizing over~\eqref{eqn:k-miblp} for \say{small} values of $k$ already yields much stronger 
dual bounds for~\eqref{eqn:miblp-vf}. In this work, the relevance of this hierarchy of relaxations is mainly theoretical. 
We show in Section~\ref{sec:valid} that the feasible solutions of~\eqref{eqn:k-miblp} can be described
by adding to the MILP relaxation $\S$ a specific (finite) class of valid linear inequalities arising from improving directions 
(or equivalently, improving solutions).
From a theoretical standpoint, this gives us some insight regarding the strength
of such a class of inequalities.

\begin{example} 
Figure~\ref{fig:3dex_k_opt} shows slices of $\P$ of the three-dimensional example for the 
four integer values of $x$. For each slice, we illustrate the feasible points of~\eqref{eqn:k-miblp}
for different values of $k$. In particular, we can observe in Figure~\ref{fig:3dex_x_1} 
($\hat{x} = 1)$ that
\begin{align*}
            \RR(\hat{x}, 1) \setminus \RR(\hat{x}) &= \brak{(\hat{x}, 3, 2), (\hat{x}, 7, 3), (\hat{x}, 2, 2), (\hat{x}, 1, 2)} \supset\\
    \supset \RR(\hat{x}, 2) \setminus \RR(\hat{x}) &= \brak{(\hat{x}, 2, 2), (\hat{x}, 1, 2)} \supset \\
    \supset \RR(\hat{x}, 3) \setminus \RR(\hat{x}) &= \brak{(\hat{x}, 1, 2)}, \\
\end{align*}
and that
$$ \F(1) \supset \F(2) \supset \F(3) \supset \F(4) = \F \cup \brak{(3, 4, 1)}.$$
Note that $(3, 4, 1)$ is the only element in $\F(4)$ that is not also in $\F$, since for its unique IFD $w = (4, -1)$ we have $\norm{w}_1 = 5$. 
Therefore, for all $k \geq 5$, $\F(k) = \F$.

\begin{figure}[h]
    \centering
    \begin{subfigure}[b]{0.45\textwidth}
        \centering
        \includegraphics[width=\linewidth]{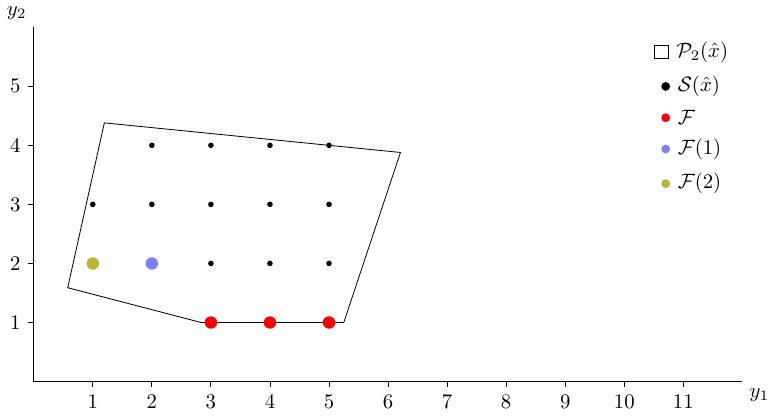} 
        \caption{$\hat{x} = 0$}
        \label{fig:3dex_x_0}
    \end{subfigure}
    \hfill
    \begin{subfigure}[b]{0.45\textwidth}
        \centering
        \includegraphics[width=\linewidth]{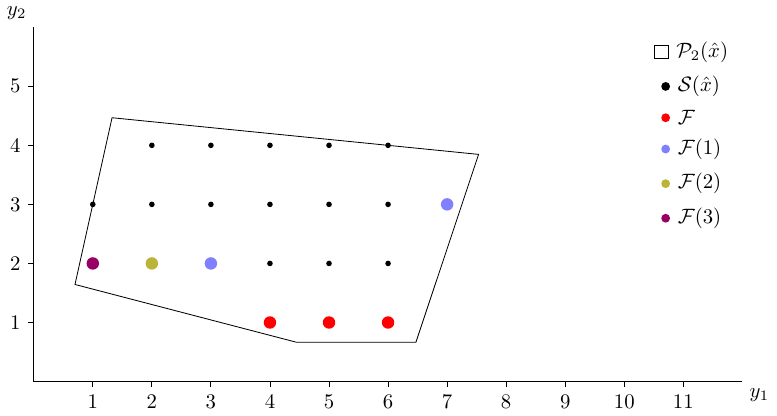} 
        \caption{$\hat{x} = 1$}
        \label{fig:3dex_x_1}
    \end{subfigure}

    \begin{subfigure}[b]{0.45\textwidth}
        \centering
        \includegraphics[width=\linewidth]{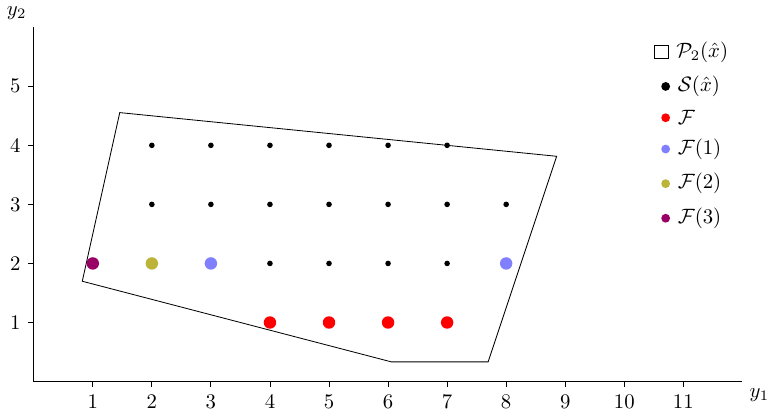} 
        \caption{$\hat{x} = 2$}
        \label{fig:3dex_x_2}
    \end{subfigure}
    \hfill
    \begin{subfigure}[b]{0.45\textwidth}
        \centering
        \includegraphics[width=\linewidth]{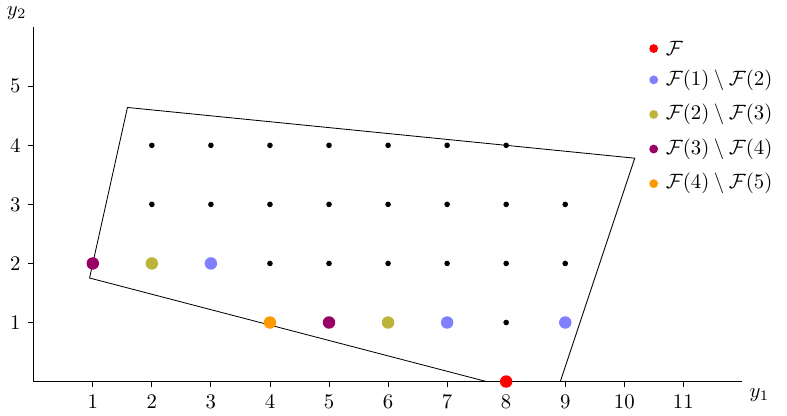} 
        \caption{$\hat{x} = 3$}
        \label{fig:3dex_x_3}
    \end{subfigure}
    \caption{Slices of the feasible region of~\eqref{eqn:k-miblp}.}
    \label{fig:3dex_k_opt}
\end{figure}

\end{example}

\section{Valid Inequalities from IDs}
\label{sec:valid}

In this section, we discuss the generation of valid inequalities from
improving directions. Our goal is to formalize what seems to be an intuitive
connection between IDICs and the $k$-opt hierarchy introduced in
Section~\ref{sec:k-opt}, both of which are derived from the concept of
improving directions. The connection we aim to establish is that, roughly
speaking, $\F(k)$ can be described using IDICs derived only from BFSs obtained
from improving directions $w \in \W$ and such that $\norm{w}_1 \leq k$. After
briefly reviewing the definitions related to ICs in Section~\ref{sec:ICs} and
introducing a notion of rank analogous to that for valid inequalities in MILPs
in Section~\ref{sec:IC-closure}, we prove the main theoretical result in
Section~\ref{sec:k-IDICs}.

\subsection{Intersection Cuts \label{sec:ICs}}

As already mentioned, we focus on two related types of intersection cuts that
can be derived from improving directions and improving solutions,
respectively. Specifically, the classes differ in the definition of bilevel free set.

\paragraph{Improving Direction Intersection Cuts.}

We first describe the BFS used for generating IDICs.
\begin{theorem}[\cite{FisLjuMonSinUse18}]\label{thm:idic}
    Let $(\hat{x},\hat{y}) \in \Re^{n_1 + n_2}$ be the extreme point of a
    simplicial radial cone $\R(\hat{x}, \hat{y})$ containing
    $\F \not= \emptyset$ and $w\in \W(\hat{x},\hat{y})$. Then we have
    that \begin{equation*} \alpha^x x + \alpha^y y \geq \beta \quad \forall
    (x, y) \in \F, \end{equation*} where this inequality is the IC generated
    from the bilevel free
    set \begin{equation}\label{eqn:idicset}\tag{IDIC} \Cid{w}
    = \left\{(x,y)\in\Re^{n_1} \times \Re^{n_2}\midd A^2 x+ G^2(y+w)\geq b^2 -
    1, y+w\geq -1\right\}. \end{equation} Furthermore, $\alpha^x \hat{x}
    + \alpha^y \hat{y} < \beta$.
\end{theorem}
To see why $\Cid{w}$ is a BFS for any $w \in \W$, observe that for
$(\hat{x}, \hat{y}) \in \S \cap \operatorname{int}(\Cid{w})$, $\hat{y} + w$ is
feasible for the associated follower's problem, which means that $w$ is an IFD
with respect to $(\hat{x}, \hat{y})$. Hence, $(\hat{x}, \hat{y})$ must be
bilevel infeasible. Thus, the set $\Cid{w}$ contains all points with respect
to which $w$ is an improving feasible direction.

Note that the definition of BFS given in~\eqref{eqn:idicset} is independent of
the point $(\hat{x}, \hat{y})$. As long as $w$ is an ID, $\Cid{w}$ is a
BFS---it need not be feasible with respect to any particular point (although
the BFS associated with a given ID could be empty). The reason we may want to
construct a direction that is an IFD w.r.t. a specific point is that this
ensures the point lies in the interior of $\Cid{w}$, which in turn ensures
that the point will violate the generated IC.

To summarize, separation of $(\hat{x}, \hat{y}) \in \P \setminus \F$ by some IC can
be guaranteed if 
\begin{itemize}
    \item[$(i)$] $\W(\hat{x}, \hat{y}) \not= \emptyset$;
    \item[$(ii)$] we can construct a simplicial radial cone
    $\R(\hat{x}, \hat{y}) \supseteq \F$ with $(\hat{x}, \hat{y})$ as its
    extreme point; and
    \item[$(iii)$] $\R(\hat{x}, \hat{y}) \not\subseteq \Cid{w}$ ($\F \not= \emptyset$).
\end{itemize}
By Proposition~\ref{lemma:bfeas}, condition $(i)$ is automatically satisfied
whenever $(\hat{x}, \hat{y}) \in \S \setminus \F$. In practical computation,
the case that $\W(\hat{x}, \hat{y}) = \emptyset$ may arise and is an important
consideration, as discussed further in~\cite{TahRal25}. In such a case, we
cannot separate $(\hat{x}, \hat{y}) \in \P$ from $\F$ with an IDIC, although
we can still do so with an inequality valid for $\conv(\S)$. Condition
$(ii)$ is typically easy to satisfy, since when $(\hat{x}, \hat{y})$ is an
extreme point of $\P$, the simplicial cone arises naturally from an associated
LP basis. Violation of condition $(iii)$ means $\F = \emptyset$, which
typically only happens after branching constraints have been applied in the
context of a branch-and-cut algorithm.

\paragraph{Improving Solution Intersection Cuts.}

Let us now consider, in constrast, the BFSs used to generate ISICs.
\begin{theorem}[\cite{FisLjuMonSinNew17,FisLjuMonSinUse18}]\label{thm:isic}
    Let $(\hat{x},\hat{y}) \in \Re^{n_1 + n_2}$ be the extreme point of a
    simplicial radial cone $\R(\hat{x}, \hat{y})$ containing $\F$ such that
    $d^2\hat{y} > d^2 y^*$ for some 
    $y^* \in \P_2(\hat{x}) \cap Y$. Then, under the stated assumptions, we have
    \begin{equation*}
    \alpha^x x + \alpha^y y \geq \beta \quad \forall (x, y) \in \F,
    \end{equation*}
    where this inequality is the IC
    associated with the bilevel free set
    \begin{equation}\label{eqn:isicset} \tag{ISIC}
    \Cis{y^*} = \left\{(x,y)\in\Re^{n_1} \times \Re^{n_2}\midd 
    d^2y \geq d^2 y^*, A^2 x \geq b^2 - G^2y^* - 1 \right\}, 
    \end{equation}
    Furthermore, $\alpha^x \hat{x} + \alpha^y \hat{y} < \beta$. 
\end{theorem}
The convex set~\eqref{eqn:isicset} includes all points $(x, y)$ for which $x$
satisfies the follower's constraints for the \emph{fixed} improving solution
$y^* \in Y$ and for which $y$ has a second-level objective no better than that
of $y^*$. As with IDICs, the BFS does not depend on the point
$(\hat{x}, \hat{y})$. There is a BFS associated with each $y^* \in Y$ (though
again, some could be empty). The reason we may desire a $y^*$ such that $d^2
y^* < d^2\hat{y}$ is to guarantee that $(\hat{x}, \hat{y})$ can be separated.

As with IDICs, it may be possible to generate an inequality when
$(\hat{x}, \hat{y}) \not\in \S$. Conditions $(ii)$ and $(iii)$ for
separation by an IDIC must also be satisfied for separation by an ISIC, but
instead of condition $(i)$, we must have $\P_2(\hat{x}) \cap
Y \not= \emptyset$. By Proposition~\ref{lemma:bfeas}, this is assured when
$(\hat{x},\hat{y}) \in \S \setminus \F$.

\begin{example} 
Let us consider the example from~\citt{MooBarMixed90}. The optimal solution
$(\hat{x}, \hat{y})$ to the LP relaxation satisfies integrality requirements,
but is bilevel infeasible. In this case, it is straightforward to see that $w
= -1$ is an improving feasible direction and $y^* = 2$ is an improving
solution. Figure~\ref{fig:ic} depicts all possible IDICs and ISICs obtained by
combining the direction and solution with both BFSs. For the sake of this
example, the BFSs reported here are defined by the original constraints rather
than the relaxed right-hand sides. Although separation of the current solution
is guaranteed under broad assumptions, for deeper cuts, we want $\Cid{w}$ to
be as large as possible and this means choosing a \say{short} directions $w$.
On the other hand, larger set $\Cis{y^*}$ arise from solutions $y^* \in Y$ that
are \say{further} from $\hat{y}$, i.e., for which $\norm{y^* - \hat{y}}$
is \say{larger}.

\begin{figure}[!ht]
    \centering
    \begin{subfigure}{0.49\textwidth}
        \includegraphics[width=\linewidth]{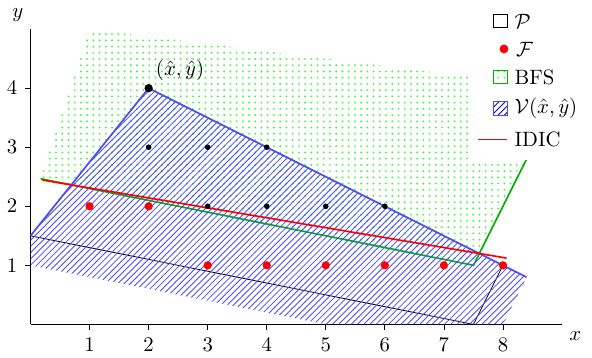}
        \caption{IDIC with IFD $w = -1$}
        \label{fig:idic:1}
    \end{subfigure}    
    \begin{subfigure}{0.49\textwidth}
        \includegraphics[width=\linewidth]{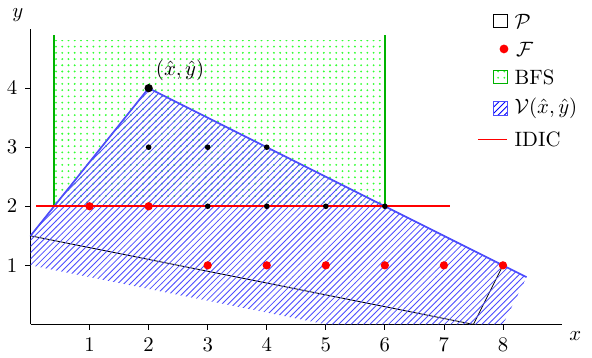}
        \caption{ISIC with IS $y^* = 2$}
        \label{fig:isic:1}
    \end{subfigure}  
    \begin{subfigure}{0.49\textwidth}
        \includegraphics[width=\linewidth]{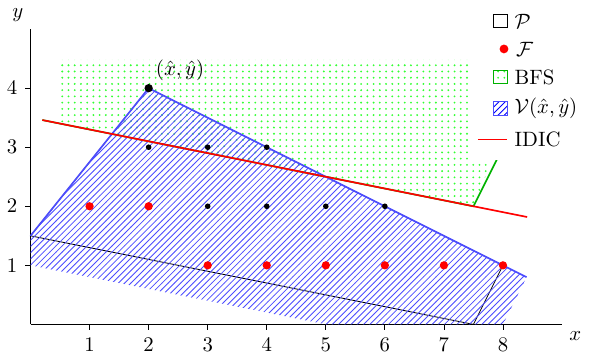}
        \caption{IDIC with IFD $y^* - \hat{y}$}
        \label{fig:idic:2}
    \end{subfigure}    
    \begin{subfigure}{0.49\textwidth}
        \includegraphics[width=\linewidth]{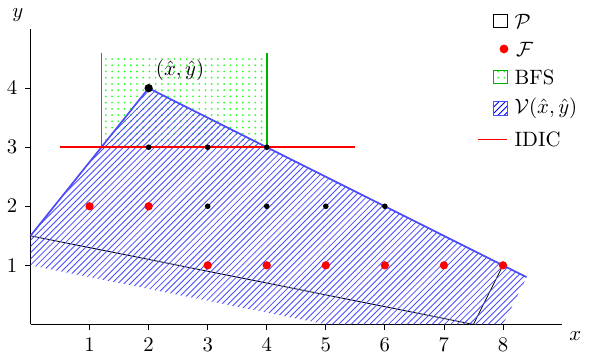}
        \caption{ISIC with IS $\hat{y} + w$}
        \label{fig:isic:2}
    \end{subfigure}    
    
    \caption{Illustration of IS and IDICs using different improving solutions/directions
    on Moore and Bard example~\citp{MooBarMixed90}}
    \label{fig:ic}
\end{figure}
\end{example}

\paragraph{Connections between ISICs and IDICs.} 
For any given $(\hat{x},\hat{y}) \in \S \setminus \F$, it is easy to verify
that if we are given any improving direction $w \in \W(\hat{x},\hat{y})$, we
can use it to obtain an improving solution $y^* = \hat{y} +
w \in \P_2(\hat{x}) \cap Y$ and vice versa. However, this equivalence does not
always hold for points in $\P \setminus \F$. Consider, for instance, the point $(1,
2.2) \in \P \setminus \S$ from the~\citt{MooBarMixed90} example.
Figure~\ref{fig:idic:1} shows that $(1, 2.2) \notin \Cid{-1}$, implying $\W(1,
2.2) = \emptyset$, even though Figure~\ref{fig:idic:2} confirms that $(1,
2.2) \in \Cis{2}$. Thus, there may exist points in $\P \setminus \S$ that
admit improving solutions, but no improving feasible directions. In this case,
it occurs because the direction $w = -0.2$ does not satisfy the integrality
requirements. As a consequence, the separation of fractional extreme points
may fail even when $d^2\hat{y} > \phi(b^2 - A^2\hat{x})$, whereas this cannot
happen with ISICs.

Another important aspect is the contrast of the two BFSs. The set $\Cid{w}$ is
defined by a direction $w$, not a fixed second-level solution. In contrast,
$\Cis{y^*}$ is defined with respect to a fixed $y^* \in Y$. Intuitively, to
enlarge $\Cis{y^*}$, $y^*$ should be a \emph{high quality} solution to the
second-level problem---corresponding to directions of larger magnitude---which
is in opposition to the goal of choosing $w$ such that $\Cid{w}$ is large. For
this reason, depending on what kind of IC we are generating, different
directions must be considered desirable. This issue is addressed in depth in
Section~\ref{sec:genid}.

\subsection{Closures and Rank for Intersection Cuts \label{sec:IC-closure}} 

In the context of MILPs, one way of characterizing the overall strength of a
specific class of valid inequalities is by analyzing the
so-called \emph{closure}, which is the convex set obtained by adding all
inequalities in the class to the initial LP relaxation. Before making the
connections between IDICs and the $k$-opt hierarchy more formal, we first show
how to apply the standard notions of closure and rank to the ICs discussed in
Section~\ref{sec:ICs}.

As with MILPs, when the set of non-dominated inequalities in a class is
finite, the closure is a polyhedron and is itself a relaxation of the original
problem. The inequalities defining this first closure are defined to have rank
1. Taking the closure again with respect to the relaxation defined by the
first closure yields the rank 2 closure and this process can be iterated. In
general, the closure of rank $r$ is the closure with respect to that of rank
$r-1$. The rank of a given valid inequality with respect to this hierarchy is
the smallest value of $r$ such that the inequality is valid for the rank $r$
closure, but not the rank $r-1$ closure~\citp{CorValid08}.
 
To apply these concepts to ICs, an obvious approach would be to consider all
cuts that can be derived from the procedure of Definition~\ref{def:IC}, taking
$(\hat{x}, \hat{y})$ to be one of the extreme points of $\P$. Given an extreme
point $(\hat{x}, \hat{y})$, the definition requires specifying a simplicial
radial cone pointed at $(\hat{x}, \hat{y})$. Such a cone is easily obtained
from a basis of the LP relaxation with respect to which $(\hat{x}, \hat{y})$
is the associated basic feasible solution. Different bases yield different
cuts, so a closure could be derived by considering all cones arising from all
bases for all extreme points of $\P$ and combining these with all possible
BFSs.

A simpler construction is one in which we enumerate all possible BFSs and
consider the convex hull of the complement of the interior for each. 
By taking the intersection of all such complements with $\P$, we obtain a similar
closure. More formally, for each $w \in \W$, let $D^1(w) :
= \conv(\P \setminus \inter(\Cid{w}))$. Then the first closure of $\P$ with
respect to the class of ICs associated with BFSs~\eqref{eqn:idicset} is called
the \emph{rank 1 IDIC closure} and defined as
$$ \P^1_{\rm ID} = \bigcap_{w \in \W} D^1(w).$$ 

Although the set of all IDs $\W$ is not finite in general, it can be replaced
by a finite subset consisting of IDs with 1-norm at most $\bar{k}$ (as defined
earlier) in the above. The \emph{rank $r$ IDIC closure}, denoted by
$\P^r_{\rm ID}$, is recursively defined as the rank 1 closure of
$\P^{r-1}_{\rm ID}$ as follows
$$ \P^r_{\rm ID} = \bigcap_{w \in \W} D^r(w), $$
where $D^r(w) := \conv(\P^{r - 1}_{\rm ID} \setminus \inter(\Cid{w}))$,
for all $w \in \W$, and $\P^0_{\rm ID} := \P$.
Hence, the closure yields a natural hierarchy of relaxations over these classes of valid
inequalities, ranging from low-rank (weaker) to high-rank (stronger) cuts.

As noted before, since there are points in $\P \setminus \conv(\F)$ that
cannot be separated by IDICs at all, the facet-defining inequalities of
$\conv(\F)$ cannot all be expected to have finite IDIC rank. The next example
shows a case in which there exists an $r^*$ such that $\P^{r^*+1}_{\rm ID}
= \P^{r^*}_{\rm ID}$ and $\P^{r^*}_{\rm ID} \cap \S \not= \F$.

\begin{example}
    Figure~\ref{fig:frac_polytope} depicts a polytope $\P \supset \F$ of
    the~\citt{MooBarMixed90} example. As illustrated, all extreme points $(x,
    y) \in \operatorname{ext}(\P)$ are fractional, i.e., $(x, y) \notin \S$.
    However, none of these extreme points lie in $\inter(\Cid{-1})$ (depicted
    in green), and thus cannot be separated by an IDIC. Consequently, we have
    that $\P^r_{\rm ID} = \P$, for all $r > 0$. Moreover, $\P^r_{\rm
    ID} \cap \S \not= \F$, as the point $(3, 2) \in (\P^r_{\rm
    ID} \cap \S) \setminus \F$.
    
    \begin{figure}[!ht]
        \centering
        \includegraphics[width=0.7\linewidth]{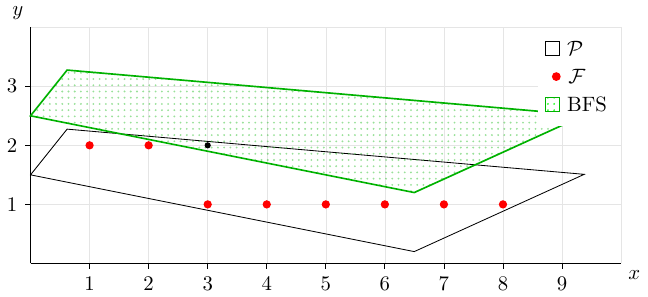} 
        \caption{Illustration of a polytope where all extreme points are fractional 
        and have an empty set of IFDs.}
        \label{fig:frac_polytope}
    \end{figure}
\end{example}

On the other hand, it is possible to separate any point in $\S \setminus \F$
that is an extreme point of $\P$ from $\conv(\F)$ with an IDIC, so a pure
cutting plane method using a combination of IDICs and MILP cuts is possible.
Next, we show that $\P^r_{\rm ID}$ is a polyhedron for any $r > 0$.
\begin{proposition}
    \label{prop:polyhedron}
    For all $w \in \W$, $D^1(w)$ is a polyhedron.
\end{proposition}

\myproof{We show that the set of extreme points of $\conv(\P \setminus \inter(\Cid{w}))$ 
is finite. Note that $\P \cap \Cid{w}$ is a polyhedron. 
Hence, $\operatorname{ext}(\P \cap \Cid{w})$ is a finite set. 
Then, by construction we have that
$$ \operatorname{ext}(\conv(\P \setminus \inter(\Cid{w}))) \subseteq 
\operatorname{ext}(\P) \cup \operatorname{ext}(\P \cap \Cid{w}).$$
Since both sets on the right-hand side are finite, so must be the left-hand side.}

\begin{theorem}
    \label{thm:closure}
    Let $r > 0$ be given. Then $\P^r_{\rm ID}$ is a polyhedron.
\end{theorem}

\myproof{ 
Note that for all $\hat{w} \in \W \setminus \W^{\bar{k}}$, we have 
$\conv(\P \setminus \inter(\Cid{\hat{w}})) = \P$. 
Then 
$$ \P^r_{\rm ID} = \bigcap_{w \in \W} D^r(w) = 
\bigcap_{w \in \W^{\bar{k}}} D^r(w),$$
for all $r > 0$. First, we prove the statement for $r = 1$. For all
$w \in \W^{\bar{k}}$, $D^1(w)$ is polyhedral by Proposition~\ref{prop:polyhedron}. 
It follows that $\P^1_{\rm ID}$ is obtained as the intersection of a finite
number of polyhedra. Assuming the statement holds for $r - 1$, we prove the
statement for $r > 1$. 
Since $\P^{r - 1}_{\rm ID} \subseteq \P$ is polyhedral, by Proposition~\ref{prop:polyhedron}
we have that $D^r(w)$ is still polyhedral, for all $w \in \W^{\bar{k}}$.
Hence, $\P^r_{\rm ID}$ is the result of the intersection of a finite number
of polyhedra.}

The same analysis can be applied to the ISICs defined by~\eqref{eqn:isicset}
to obtain a hierarchy defined by closures $\P^r_{\rm IS}$.

\subsection{Intersection Cuts and $k$-optimality \label{sec:k-IDICs}}

The hierarchy of relaxations presented in the previous section provides a
framework to classify the strength of IDICs as a function of their IDIC rank,
which is the natural counterpart of the similar hierarchies in the theory of
valid inequalities for MILPs (Chv\'atal-Gomory rank, etc.). Our goal in this
section, however, is to analyze a different and more practical hierarchy 
that classifies the strength of IDICs as a function of the points in $\S$
they can separate, with a close relationship to the $k$-opt hierarchy 
introduced earlier.

The analysis in this section was the original motivation for undertaking the 
work presented in this paper, driven by the desire of computing
the strong bounds produced by the $k$-opt relaxation.
As our main result in this section shows, these bounds can be obtained with only 
a minimal modification to the existing branch-and-cut framework of~\citt{TahRalDeN20} 
by restricting the generation of valid inequalities to the subset of IDICs that are 
valid for the feasible points of the $k$-opt hierarchy, that is, those derived
from IFDs with 1-norm at most $k$.
To formalize this result, we first introduce the notion of a \emph{$k$-IDIC},
an IDIC valid for~\eqref{eqn:k-miblp}.
\begin{definition}
    \label{def:kIC}
    For $k \in [\bar{k}]$, a $k$-IDIC is an IC generated from a BFS $\Cid{w}$
    such that $w \in \W^k$. 
\end{definition}
Before getting into the formalities, we first outline the intuition behind the
proof. Consider any bilevel infeasible point $(x, y) \in \S \setminus \F$ and
let $k^* = \min_{w \in \W(x, y)} \norm{w}_1$ be the smallest 1-norm of any IFD
w.r.t. $(x, y)$. Proposition~\ref{lemma:k-bfeas} and
Theorem~\ref{thm:hierarchy} together imply that $(x, y) \in \F(k)$ if and only
if $k \leq k^* - 1$. In particular, we must have $\W(x, y;
k^*) \not= \emptyset$ and hence, there must be at least one IFD $w \in \W(x,
y; k^*)$. From this IFD $w$, given an appropriate simplicial cone (which we
show below always exists), we can always derive a $k^*$-IDIC separating $(x,
y)$ from $\F(k^*)$ (and hence also from $\F(k)$ for all $k^* \leq
k \leq \bar{k}$). 

To state the same logic in other terms, Proposition~\ref{lemma:k-bfeas} and
Theorem~\ref{thm:hierarchy} together specify a partition of $\S \setminus \F$
into \emph{levels} that correspond to levels of the $k$-opt hierarchy. A point
$(x, y)$ is on \emph{level} $k$ if it is in $\F(k) \setminus \F(k + 1)$. 
A point on level $k$ can be separated from $\F(k + 1)$ by a
$(k+1)$-IDIC. 

It is important to emphasize that there is no relationship in general between
the IDIC \emph{rank} of an inequality (as defined in
Section~\ref{sec:IC-closure}) and the level of the points it can
separate. In fact, for any fixed $k$, a $k$-IDIC can have any IDIC rank. There
exists an entire rank hierarchy for $k$-IDICs for any fixed $k$ that mirrors
that derived in the previous section. Applying Definition~\ref{def:kIC} to the
extreme points of $\P$, for example yields $k$-IDICs of ``rank 1'' and we can
iterate to derive $k$-IDICs of higher rank.

Next, we state the main result of this section which comprises two parts. The
first establishes that $k$-IDICs are valid for $\F(k)$ by showing that for all
$w \in \W^k$, $\Cid{w}$ is ``$\F(k)$-free''.
The second shows that any point in $\S$ that is not in $\F(k)$ can be
separated by a $k$-IDIC.

\begin{theorem}
    \label{thm:k-separation}
    For all $k \in [\bar{k}]$, the following statements hold:
    \begin{itemize}
        \item[$(1)$] all $k$-IDICs are valid for $\F(k)$; and
        \item[$(2)$] for all $(x, y) \in \S \setminus \F(k)$, there exists a $k$-IDIC violated by $(x, y)$.
    \end{itemize}
    In other words, $k$-IDICs can separate \emph{all and only} points in $\S \setminus \F(k)$.
\end{theorem}

\myproof{
$(1)$ Let $k \in [\bar{k}]$ and $w \in \W^k$ be
given. We show that $\inter(\Cid{w}) \cap \F(k) = \emptyset$.
By definition of~\eqref{eqn:idicset}, $w \in \W(x, y; k)$ for all $(x, y) \in \inter(\Cid{w}) \cap \S$.
Consequently, Proposition~\ref{lemma:k-bfeas} implies $(x, y) \notin \F(k)$.

$(2)$ Let $k \in [\bar{k}]$ and $(x, y) \in \S \setminus \F(k)$ be given.
First, we show that a suitable simplicial radial cone can be constructed 
from standard geometric arguments.
The polyhedron 
$$\F^{(x, y)}(k) = \conv(\F(k) \cup \{(x, y)\})$$
contains $\F(k)$ and has $(x, y)$ as one of its extreme points. Moreover, 
there exists a simplicial radial cone $\R(x, y)$ pointed at $(x, y)$ 
whose description is obtained by selecting a
set of linearly independent inequalities binding at $(x, y)$ from the
description of $\F^{(x, y)}(k)$.
Then, by Proposition~\ref{lemma:k-bfeas} there exists a $w \in \W(x, y; k)$;
and $(x, y) \in \inter(\Cid{w})$ by the definition of~\eqref{eqn:idicset}.
Therefore,
$$\conv\left(\F^{(x, y)}(k) \setminus \inter(\Cid{w})\right)$$
is a $k$-IDIC violated by $(x, y)$.
}

To complete the argument, we show that, for any given $k \in [\bar{k}]$, 
the convexification of the set of points satisfying both integrality 
requirements and all possible $k$-IDICs derived from Theorem~\ref{thm:k-separation}
coincides with that of the feasible region of~\eqref{eqn:k-miblp}.
As such, for all $k \in [\bar{k}]$, let
$$\Pi^k_{\rm ID} = \bigcap_{(x, y) \in \F(k - 1) \setminus \F(k)} \bigcap_{w \in \W(x, y; k)} \conv\left( \F^{(x, y)}(k) \setminus \inter(\Cid{w})\right)  .$$
be the intersection of all the $k$-IDICs separating from $\F(k)$ points in $\F(k - 1) \setminus \F(k)$, those having level $k - 1$.

\begin{corollary}
    \label{cor:convexify}
    For all $k \in [\bar{k}]$, $\conv(\S \cap \Pi^k_{\rm ID}) = \conv(\F(k))$.
\end{corollary}
    
From a theoretical standpoint, the previous result connects the strength of
IDICs generated from an ID with 1-norm $k$ to the $k$-opt relaxation. From a
practical standpoint, it suggests that the separation of this family of cuts
produces dual bounds converging to those associated with the
hierarchy~\eqref{eqn:k-miblp}. As a matter of fact, branch-and-cut and
cutting-plane methods can be interpreted as optimizing over the convex hull of
a certain feasible region. Therefore, such dual bounds can, in principle, be
computed by a pure cutting plane in which integrality is restored with the
separation of standard MILP cuts, while $k$-optimality is enforced by
$k$-IDICs. Equivalently, the same bounds can be produced by a branch-and-cut
algorithm exclusively generating $k$-IDICs, with integrality enforced
(convexification) through standard MILP branching.

We end the section with a few complementary results, illustrated through examples.
In contrast with the properties established for IDICs, an analogous result does not
hold in general for ISICs. Specifically, given a point
$(\hat{x}, \hat{y}) \in \S \setminus \conv(\F(k))$ and a
$w \in \W(\hat{x}, \hat{y}; k)$, it is not necessarily the case that $w$ is an IFD for
every $(x, y) \in \inter(\Cis{\hat{y} + w}) \cap \S$, as we exemplify next.

\begin{figure}[!ht]
    \centering
    \includegraphics[width=0.7\linewidth]{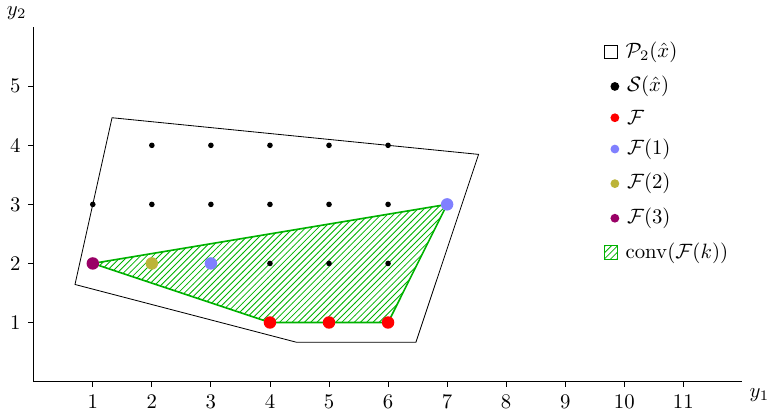}
    \caption{A slice of $\conv(\F(1))$ of the three-dimensional example, for $\hat{x} = 1$}
    \label{fig:conv-k}
\end{figure}

\begin{example}
    Figure~\ref{fig:conv-k} shows the slice of the three-dimensional example for $\hat{x} = 1$. 
    Let us consider $\hat{y} = (4, 3)$, $w = (0, -1)$ and $y^* = (4, 2)$ as the improving solution.
    The BFS is defined as $\Cis{y^*} = \brak{(x, y) \in \Re^{n_1} \times \Re^{n_2} \midd y_2 \geq 2, 0 \leq x \leq 3}$.
    One can verify that $(1, 7, 3) \in \Cis{y^*}$, but $w \notin \W(1, 7, 3; 1) = \emptyset$.
\end{example}

Note that, unlike the MILP setting where all points satisfying integrality 
requirements within the convex hull of $\S$ are feasible, there may 
exist points $(x, y) \in \conv(\F(k)) \cap \S$ for which $\W(x, y; k) \not = \emptyset$, 
and thus $(x, y) \notin \F(k)$, as illustrated in the following example.

\begin{example}
    Figure~\ref{fig:conv-k} shows a slice of $\conv(\F(1))$ (in green) of the 
    three-dimensional example, for $\hat{x} = 1$. It is easy to verify that the vector 
    $\hat{w} = (0, -1)$ is an IFD for the point $(4, 2) \in \S(\hat{x})$, i.e., $\W(4, 2; k) \not = \emptyset$.
    Although $(4, 2) \in \inter(\Cid{\hat{w}})$, there exists no simplicial radial cone pointed at 
    $(4, 2)$ that contains $\F(k)$, since $(4, 2) \in \conv(\F(k))$. 
    Therefore, in this case, separation with $k$-IDICs is not possible but also unnecessary.
\end{example}

\section{Branch-and-Cut Algorithm}
\label{sec:algo}

In this section, we present a branch-and-cut algorithm based on the
overarching framework of~\citp{TahRalDeN20,TahRal25} but in which the
feasibility check and the generation of valid inequalities are unified with an
oracle for finding an IFD. The implementation is based on the open source
solver \MIBS{}~\citp{MIBS1.2} and we refer the reader
to~\citp{TahRalDeN20,TahRal25} for details. This section describes the details
of the main components, focusing on the differences from what is currently
implemented in \MIBS{} \t{1.2}.

\subsection{General Framework}

The outline of the overall method is in Algorithm~\ref{algo:scheme}.
\begin{algorithm}[ht]
    \caption{Generic Branch-and-Cut Using Improving Direction Oracle}
    \label{algo:scheme}
    \begin{algorithmic}[1]
        \State Initialize the set $\J$ of unexplored subproblems with the original problem.
        \State $U \gets \infty$
        \While{$\J \not = \emptyset$}
            \State \textbf{Select} a node $t$ from $\J$ with feasible region $\F^t$.
            \State \textbf{Bound} node $t$ to generate dual bound $L^t$ and candidate solution $(x^t, y^t)$. \label{step:bounding}
            \If{$L^t = \infty$ \textbf{or} $L^t \geq U$}
                \State \textbf{Prune} node $t$. \label{step:prune-infeas}
            \EndIf
            \State Check $\W(x^t, y^t) \stackrel{?}{=} \emptyset$ (optional if $(x^t, y^t) \not\in \S$). \label{step:id-gen}
            \If{$\W(x^t, y^t) \not = \emptyset$ and/or $(x^t, y^t) \not\in \S$}
                \State Either 
                \begin{itemize}
                \setlength{\itemindent}{3em}
                \item \textbf{separate} $(x^t, y^t)$ from $\conv(\F)$ and put $t$ back in $\J$; or 
                \item \textbf{branch} and add new subproblems $t_1, \dots, t_k$ to $\J$. 
                \end{itemize} \label{step:cut-or-branch}
            \Else 
                \State $(x^t, y^t) \in \F$ by Proposition~\ref{lemma:bfeas}. \label{step:feasible}
                \State $U^t \gets cx^t + d^1y^t$.
                \State $U \gets \min\{ U, U^t \}$.
                \State \textbf{Prune} node $t$. \label{step:prune-feas}
            \EndIf
             
        \EndWhile
    \end{algorithmic}
\end{algorithm}
Obviously, this is only a general framework and the precise way in which
various steps in this algorithm are implemented determines crucially how the
algorithm will perform in practice. In particular, the choice between
separation (generating a violated valid inequality) and branching on
line~\ref{step:cut-or-branch} is important. Most of the detailed control
mechanisms, including the one for making the decision between branching and
cutting, are the same here as in \MIBS{} \t{1.2} and are discussed
in~\citp{TahRal25}. In the remainder of this section, we focus on the method
used for generating IFDs, which is the main innovation.

\subsection{Generating Feasible Improving Directions}
\label{sec:genid}

Line~\ref{step:id-gen} of the algorithm is the step that differentiates our
approach from that of~\citt{TahRalDeN20}.
We now switch gears and discuss the empirical issues surrounding the
framework that has been introduced in the previous sections. 
In particular, we focus on 
the practical details of the problem of generating an appropriate IFD, given a
candidate solution $(\hat{x}, \hat{y}) \in \P$. As we already mentioned, the
goal of finding an IFD is twofold. First, the IFD (or the proof that none
exists), certifies the feasibility status of $(\hat{x}, \hat{y})$ by
Proposition~\ref{lemma:bfeas} when $(\hat{x}, \hat{y}) \in \S$ (otherwise,
$(\hat{x}, \hat{y})$ is trivially infeasible). Second (and perhaps more
crucially), the direction serves to produce an IC violated by
$(\hat{x}, \hat{y})$. For this latter purpose, the particular direction
produced makes a big difference. In the following two sections, we first
discuss relevant objective functions and then both an exact and a heuristic
way of generating directions according to these objective functions.

\subsubsection{Objective Functions}

Although any IFD in $\W(\hat{x},\hat{y})$ suffices to generate a BFS of either
the form~\eqref{eqn:idicset} or the form~\eqref{eqn:isicset} that can be used
to generate an associated IC violated by $(\hat{x}, \hat{y})$, the particular
objective function used serves to guide the selection, ideally leading to BFSs
that are \say{larger} and therefore produce deeper cuts.

Let us consider the BFS~\eqref{eqn:idicset} first. As previously remarked, a
direction $w \in \W(\hat{x},\hat{y})$ is likely to include more points in
$\Cid{w}$ if its 1-norm is small. Clearly, a possible approach is to
consider the following objective function
\begin{equation}
    \label{eqn:magnitude}
    \min \; \norm{w}_1.
\end{equation}
\citt{FisLjuMonSinUse18} proposed an objective function
that has a related goal
\begin{equation}
    \min \left\{\sum_{i=1}^{m_2} \max\{g^2_i w, 0\} + \norm{w}_1\right\}, \label{eqn:idic-obj}
\end{equation}
where $g^2_i$ is the $i^{\textrm th}$ row of $G^2$. This objective can be linearized in the obvious way and
the idea is that row $i$ can be dropped in the definition of set~\eqref{eqn:idicset} 
if $g^2_i w \leq 0$, enlarging the BFS.

For the BFS~\eqref{eqn:isicset}, although it is constructed based on a
solution, we can nevertheless consider any $w \in \W(\hat{x},\hat{y})$ and set
$y^* = \hat{y} +w$. Intuitively, a good direction for this IC is obtained by
maximizing the improvement in the follower's objective function by using the
objective 
\begin{equation}
    \min d^2w.
\end{equation}
Other objectives are also possible. For example, an objective with the
philosophy similar to that of~\eqref{eqn:idic-obj} is suggested
in~\cite{FisLjuMonSinNew17}.

\subsubsection{Algorithms \label{sec:id-gen}}

\paragraph{Exact.} Given a desired objective function, the most
straightforward way of generating an IFD is to take the approach described
in~\citp{FisLjuMonSinUse18, TahRal25} for the generation of IDICs. That is, we
describe elements of $\W(\hat{x}, \hat{y})$ as the points in the following
mixed-integer set and directly optimize over this set with the given objective
using an off-the-shelf method for solving MILPs.
\begin{align} 
    \quad d^2 w & \leq -1 \notag\\
    \quad G^2 w & \geq b^2 - A^2 \hat{x} - G^2 \hat{y} \label{eqn:feasid} \tag{ID}\\
    \quad w & \geq -\hat{y} \notag\\
    \quad w & \in Y. \notag
\end{align}
As previously noted, the infeasibility of~\eqref{eqn:feasid} certifies
emptiness of $\W(\hat{x}, \hat{y})$, while any solution is an IFD with respect
to $(\hat{x}, \hat{y})$. From a theoretical standpoint, the question of
whether there is a point satisfying~\eqref{eqn:feasid} is equivalent to that of whether 
there is an improving solution to the follower's problem,
provided that $(\hat{x}, \hat{y}) \in \S$. From a practical standpoint, there
is a difference, since \eqref{eqn:feasid} allows us to specify an objective
function that favors certain improving directions (and thus certain improving
solutions) over others. 

\paragraph{Heuristic.} Although optimizing over~\eqref{eqn:feasid} is a
straightforward approach, the main challenges in finding IFDs this way is the
computational burden associated with solving an \NPcomplexity-hard subproblem.
Fortunately, to guarantee the correctness of a hypothetical branch-and-cut
algorithm using an oracle based on Proposition~\ref{lemma:bfeas} for checking
bilevel feasibility, we only need an exact answer to the associated decision
problem when the candidate solution satisfies integrality requirements. In
cases where separation is the primary task, a heuristic approach for finding
feasible solutions of~\eqref{eqn:feasid} would suffice. Furthermore, it may be
computationally advantageous to first run heuristics, even when the solution
belongs to $\S$, resorting to an exact algorithm only if no IFD is found using
the heuristic. Recent work by~\citt{GaaLeeLjuSinTanSOCPbased24} highlights the
benefits of heuristic methods for identifying improving solutions in the
context of generating valid inequalities for integer bilevel nonlinear
problems, further supporting the promise of this approach.

We present two heuristics aimed at producing an IFD with respect to a given
$(\hat{x}, \hat{y}) \in \P$ by restricting the feasible region to, e.g.,
elements of $\W(\hat{x}, \hat{y}; k)$, for some \say{small} value of $k$. This
choice is supported by some preliminary evidence. First, the example of
Figure~\ref{fig:3dex_k_opt} suggests that $\F(k)$ is a reasonably good
approximation of $\F$ even for small values of $k$. More broadly, the results
in~\citt{XueProRal22} also indicate that $k \approx 3$ provides a good
trade-off between the quality of the dual bound of~\eqref{eqn:k-miblp} and the
computational burden required to compute it. Corollary~\ref{cor:convexify} provides a
theoretical guarantee that the effect of generating $k$-IDICs in this context
should mirror the effect of solving the $k$-opt relaxation. Perhaps most
importantly, this problem lends itself well to heuristic approaches, such as
local searches, since if a direction $w \in \W(\hat{x}, \hat{y}; k)$ exists,
then it must direct $\hat{y}$ towards points in its $k$-neighborhood.
Consequently, a promising approach is to search exclusively on directions
leading $\hat{y}$ to points in $\N_k(\hat{y})$.

The first method is described in Algorithm~\ref{alg:enum} and uses a pure
local search algorithm. Given a value of $k$, it examines all possible
directions $w$ with $ \norm{w}_1 \leq k$ and checks if they are improving and
feasible. Note that this approach allows to use any (possibly non-linear
and/or non-convex) objective function as a measure of \say{quality} of the
multiple IFDs the local search may identify.

\begin{algorithm}
    \caption{\texttt{generateNeighbors(k, (x, y), obj)}} \label{alg:enum}
    \hspace*{\algorithmicindent} \textbf{Input}: $ \mbox{\texttt{k}} \leq \bar{k},\ 
                                                  \mbox{\texttt{(x, y)}} \in \P,\ 
                                                  \mbox{\texttt{obj}} : \Z^{n_2} \rightarrow \Re$ \\
    \hspace*{\algorithmicindent} \textbf{Output}: $w^* = \argmin_{w \in \W} \mbox{\texttt{obj}}(w)$, with $\W \subseteq \W(\hat{x}, \hat{y}; k)$
    \begin{algorithmic}[1]
        \State Let $\W \gets \emptyset$ \Comment{Initialize $\W$}
        \For{$w \in \Z^{n_2}$ such that $\norm{w}_1 \leq \mbox{\texttt{k}}$} \label{for:ls}
        \If{$d^2 w  \not\leq -1$}
            \State Discard $w$ and go to line~\ref{for:ls} \Comment{$w$ is not improving}
        \EndIf
        \If{$G^2 w \not\geq b^2 - A^2 \hat{x} - G^2 \hat{y}$}
            \State Discard $w$ and go to line~\ref{for:ls} \Comment{$w$ is not feasible}
        \EndIf
        \If{$\hat{y} + w \not\geq 0$}
            \State Discard $w$ and go to line~\ref{for:ls} \Comment{$w$ is not feasible}
        \EndIf
        \State $\W \gets \W \cup \{w \}$ \Comment{$w$ is an IFD for $(\hat{x}, \hat{y})$}
        \EndFor
        \State $w^* \gets \argmin_{w \in \W} \mbox{\texttt{obj}}(w)$
        \State \textbf{return} $w^*$
    \end{algorithmic}
\end{algorithm}

The second approach is inspired by the well-known primal heuristic for
MILPs, \textit{local branching}~\citp{FisLodRepairing08}, and consists in
intersecting the feasible region of~\eqref{eqn:feasid} with the
$k$-neighborhood of $\hat{y}$, leading to the following formulation:
\begin{align} 
    \quad d^2 w & \leq -1 \notag\\
    \quad G^2 w & \geq b^2 - A^2 \hat{x} - G^2 \hat{y} \label{eqn:k-feasid} \tag{$k$-ID}\\
    \quad w & \geq -\hat{y} \notag\\
    \quad \norm{w}_1 & \leq k \notag\\
    \quad w & \in Y. \notag
\end{align}
The resulting MILP tends to be noticeably easier to solve
than~\eqref{eqn:feasid}. Moreover, if~\eqref{eqn:idic-obj} is used as
objective function, then the same artificial variables can be used to
linearize both the objective and the 1-norm constraint.

\noprint{
\subsection{Implementation}

In this section, we go into more detail on the important steps of this algorithm. 

\paragraph{Bounding.} The efficiency of Algorithm~\ref{algo:scheme} hinges
primarily on the production of strong primal and dual bounds on the optimal
solution value of each subproblem on line~\ref{step:bounding}. Due to the
branching we employ (see below), the subproblems here differ from the original
problem only in the bounds on the variables. We assume that the original
problem formulation includes upper and lower bounds on both the leader's and
follower's variables, given by vectors $u_x, l_x \in \Q^{n_1}$ and $u_y,
l_y \in \Q^{n_2}$, respectively. Then, the feasible region of the subproblem
at node $t$ is described by
\begin{equation} \label{eqn:subproblem} \tag{NODE-$t$}
    \F^t = \{(x, y) \in \F \mid l^t_x \leq x \leq u^t_x, l^t_y \leq y \leq u^t_y\},
\end{equation}
where $u^t_x, l^t_x, u^t_y, l^t_y$ are the local vectors of upper and lower
bounds arising from branching.  

At each node $t$ of the search tree, a relaxation is solved to obtain a dual bound, 
denoted by $L^t$. This relaxation can be obtained by replacing $\F$ in~\eqref{eqn:subproblem} 
with either the LP $\P$, or an MILP, whether $\S$ or
an optimality-based relaxation following the approach of~\eqref{eqn:k-miblp}, 
and then strengthening it with valid inequalities. The choice among these relaxations 
is intricate and ultimately is an empirical question, making a comprehensive discussion 
on the implications beyond the scope of this work. However, if from the one end a convex 
relaxation is typically preferred to ensure tractability, problems that are hard for 
$\Sigma_2^p$ may benefit from a more challenging non-convex, and likely stronger 
relaxation. In the case of solving an MILP, ensuring that the solution to the relaxation is a member of $\S$ means that 
Proposition~\ref{lemma:bfeas} applies at each step and will result in less branching and a 
much smaller search tree. On the other hand, part of the branching is delegated to the subsolver 
used for the MILP relaxation. Moreover, fast re-optimization is well known to be crucial
for cut generation thus making $\P$ seemingly to achieve the best trade-off. 

Primal (upper) bounds arise by producing a member of $\F^t$, 
either by heuristic methods or by showing that $(x^t, y^t) \in \F$ on line~\ref{step:feasible}. 
Whenever the solution from the relaxation at node $t$
is not bilevel feasible, we set $U^t = \infty$. This information can be aggregated
to obtain a global primal bound by $U = \min_{t \in T} U^t$,
where $T$ is the set of terminating nodes in the search tree, also referred as \emph{leafs}.

\paragraph{Pruning.} Whenever at node $t$ we have $L^t \geq U$, that node can be 
discarded (or \emph{pruned}), since the feasible region of this subproblem does not
contain any solution better than $U$. In particular this situation happens
when subproblem at node $t$ is infeasible ($L^t = \infty$) (line~\ref{step:prune-infeas}) 
or when $L^t = U^t \geq U$ (line~\ref{step:prune-feas}).

\paragraph{Cutting.} Generation of valid inequalities occurs on line~\ref{step:cut-or-branch} 
if we choose to separate rather than branching. Note that this means that we add a valid 
inequality and continue to process the same node, which is accomplished in this 
high-level algorithm by putting $t$ back in $\J$ and then selecting it again. 
The separation can be accomplished by either IDICs or ISICs.

\paragraph{Branching.} As previously mentioned, the role of branching procedure
is to create subproblems from the current node $t$. This is accomplished by
imposing a disjunction valid for $\F$, 
dividing the feasible region of the relaxation in multiple parts (generally two)
removing portions that contain no bilevel feasible solutions. In our proposed algorithm, 
we use the standard approach of branching on variables, employing one of two strategies 
provided by MibS: either prioritizing branching on the leader's variables that 
are also present in the follower's constraints (so-called \emph{linking variables}); 
or allowing branching on any variable with fractional value. The latter is the strategy 
typically used in solving MILPs, while the former is specific to MIBLPs, revolving 
around the particular role the linking variables play. The former strategy 
prioritizes branching on linking variables as long as their value are not all 
fixed by branching, even if this may mean we branch on variables with integer 
values (see~\citp{TahRalDeN20} for a detailed discussion of branching).
}

\section{Computational Results}
\label{sec:res}

In this section, we report on the extensive empirical analysis carried out to
assess the impact of the modifications to the standard branch-and-cut approach
already implemented in \MIBS{} \t{1.2}. This analysis has several related
goals. First and foremost, we want to assess the advantage of unifying the
feasibility check with the generation of valid inequalities by use of a single
oracle. Second, within the context of Algorithm~\ref{algo:scheme}, we want to
measure the computational benefit of heuristic methods for finding IFDs. The
overall objective is to determine whether the scheme presented has the
potential to improve state-of-the-art solvers.

\subsection{Implementation Details}

We build on the foundation presented in~\citt{TahRalDeN20}, which is
implemented in the open source solver \MIBS{}, distributed by the~\citt{COIComputational18}
Foundation's repository of open source projects. For
a comprehensive description of the parameters available in \MIBS{}, we refer
the reader to the documentation available
at \url{https://github.com/coin-or/MibS}. All algorithmic variations presented
have been implemented on top of version \t{1.2.1} and their employment can be
regulated using the following newly introduced parameters.
\begin{itemize}

\item \texttt{useImprovingDirectionOracle} controls whether bilevel
    feasibility is checked using the new oracle for generating improving
    directions as in Algorithm~\ref{algo:scheme}. Otherwise, \MIBS{}'s default
    oracle is used and the scheme follows to the one presented
    in~\citt{TahRalDeN20}.

\item \texttt{improvingDirectionType} controls the
    method used to find an IFD, either the local search is used, as described
    in Algorithm~\ref{alg:enum}, or the problem~\eqref{eqn:feasid} is solved as
    an MILP. Note that even when the local search is used, the solver must
    still solve~\eqref{eqn:feasid} when necessary to guarantee
    the correctness (see discussion in
    Section~\ref{sec:genid}).

\item \texttt{maxNeighborhoodSize} controls the parameter $k$ in
    Algorithm~\ref{alg:enum} (\ref{eqn:k-feasid}, resp.)
    if \texttt{improvingDirectionType} is set to 1 (0, resp.).

\item \texttt{useLocalSearchDepthLb} and \texttt{useLocalSearchDepthUb}
    regulate the use of heuristics for finding an IFD. These parameters
    restrict the execution of Algorithm~\ref{alg:enum} or \eqref{eqn:k-feasid}
    to nodes whose depth in the search tree falls within the specified lower
    and upper bounds. Otherwise, problem~\eqref{eqn:feasid} is solved.
\end{itemize}

Most of the experiments were done using the automatic default parameters
settings in MibS 1.2.1, but when \t{useImprovingDirectionOracle}
is set to \t{True}, several default behaviors are changed.
\begin{itemize}
\item The automatic setting of parameters related to cut generation is
disabled and \emph{only} the generation of IDICs is enabled. The automatic
defaults in \MIBS{} may or may not enable generation of IDICs and/or
additional classes.
\item The \t{fractional} branching strategy is always used, even for
interdiction problems. This is because we conjecture that
the \t{fractional} strategy has a clear potential to be more  
effective in the context of Algorithm~\ref{algo:scheme}, 
since it fosters the satisfaction of integrality 
requirements at both levels, allowing a more frequent application 
of Proposition~\ref{lemma:bfeas} and, in turn, the discover of primal 
bounds. 
\item Finally, in \MIBS{} \t{1.2}, the bilevel feasibility check precedes 
the cut generation and is implemented with an oracle that
evaluates~\eqref{eqn:phi} by solving the follower's problem to optimality. By
default, this check is undertaken when either the current solution satisfies
integrality at both levels or all the linking variables are fixed. In the
latter case, \MIBS{} solves an upper bounding problem (one additional MILP) to
find the best solution for the given set of linking variables.
When \t{useImprovingDirectionOracle} is \t{True}, the usual check is bypassed,
and the feasibility check is instead undertaken \textit{after} cut generation,
since Proposition~\ref{lemma:bfeas} says that the integer solution is feasible 
if and only if we are not able to find an IFD (generate an IC).
\end{itemize}

\subsection{Dataset}
The selection of a diversified collection of instances is one of the crucial
tasks for an insightful empirical analysis since special classes of problems
can show very different behavior.
Recently, \citt{ThuKleLjuRalSch25} released the Bilevel
Optimization (Benchmark) Instance Library (BOBILib) whose intent is to provide
access to the community to a large and well-curated set of test instances, in
a similar fashion of the MIPLIB~\citp{GleHenGamAchBasBerChrJarKocLin.ea21} for
MILPs. The instances used in this work are drawn from~\citt{TahRal25}, which
includes some of the instances available in BOBILib, along with additional
ones constructed by the authors. Moreover, a subset of instances
from~\citt{XueProRal22} is included. Table~\ref{tab:datasets} provides details
on each dataset, including the number of instances, types of leader's and
follower's variables, the range of variables and constraints at each level,
and the alignment of the objective functions of the two levels. There are a
total of 399 interdiction problems and 277 IBLPs. The final count of instances
is 676.

\begin{table}[]
    \centering
    \begin{tabular}{cccccccc}
    \toprule
    \textbf{Class} & \textbf{Data Set} & \textbf{\#} & \textbf{Var. Type} & \textbf{Var\#} & \textbf{Constr\#} & \textbf{Align} & \textbf{Source} \\
    \midrule
    \multirow{4}{*}{Interdiction}&\multirow{2}{*}{INT-DEN} & \multirow{2}{*}{300} & B & 10-40 & 1     & \multirow{2}{*}{-1} & \multirow{2}{*}{\citt{DeNInterdiction11}} \\
                                 &                         &                      & B & 10-40 & 11-41 &                     &  \\
    \cline{2-8}
                                &\multirow{2}{*}{INT-SHI} &  \multirow{2}{*}{99} & B & 15-30 & 1     & \multirow{2}{*}{-1} & \multirow{2}{*}{\citt{XueProRal22}} \\
                                &                         &                      & B & 15-30 & 16-31 &                     &  \\
    \midrule
    \multirow{10}{*}{IBLP}      & \multirow{2}{*}{DEN}     & \multirow{2}{*}{50}  & I & 5-15  & 0  & \multirow{2}{*}{Varies} & \multirow{2}{*}{\citt{DeNInterdiction11}}  \\
                                &                         &                      & I & 5-15  & 20 &                         & \\
    \cline{2-8}
                                &\multirow{2}{*}{DEN2}    & \multirow{2}{*}{110} & I & 5-10  & 0 & \multirow{2}{*}{Varies} & \multirow{2}{*}{\citt{DeNInterdiction11}} \\
                                &                         &     & I & 5-20  & 5-15 &                                       & \\
    \cline{2-8}
                                &\multirow{2}{*}{ZHANG}   & \multirow{2}{*}{30}  & B & 50-80  & 0   & \multirow{2}{*}{0.6-0.8} & \multirow{2}{*}{\citt{ZhaOzaBranchcut17}} \\
                                &                         &                      & I & 70-110 & 6-7 &                          & \\
    \cline{2-8}
                                &\multirow{2}{*}{ZHANG2}  & \multirow{2}{*}{30}  & I & 50-80  & 0   & \multirow{2}{*}{0.6-0.8} & \multirow{2}{*}{\citt{ZhaOzaBranchcut17}} \\
                                &                         &                      & I & 70-110 & 6-7 &                          & \\
    \cline{2-8}
                                &\multirow{2}{*}{FIS}     & \multirow{2}{*}{57}  & B & Varies & Varies & \multirow{2}{*}{-1} & \multirow{2}{*}{\citt{FisLjuMonSinUse18}} \\
                                &                         &                      & B & Varies & Varies &                     &  \\
    \bottomrule
    \end{tabular}
    \caption{Summary of the datasets}
    \label{tab:datasets}
\end{table}

\subsection{Experimental Setup}

All experiments were conducted on compute nodes running the 
Linux (Debian 8.11) operating system with dual AMD Opteron 6128 
processors and 32 GB of RAM. All experiments were run sequentially 
with a time limit of 3600 seconds and a memory limit of 16 GB. 
Unless stated otherwise, all other parameters of \MIBS{} \t{1.2.1} were 
set to their default values.

\subsection{Configurations}

We compared various configurations using the approach of
Algorithm~\ref{algo:scheme} (\texttt{useImprovingDirectionOracle} set
to \t{True}) with various configuration of \MIBS{} \t{1.2.1}. 
Variations of Algorithm~\ref{algo:scheme} are referred to by names prefixed with 
\texttt{idB\&C}. The configurations of \texttt{idB\&C} differ in the settings
for the previously introduced parameters. Because of the numerous policies
that have been tested, we group all configurations under the following
classes, where \say{\texttt{*}} stands as a placeholder for the actual value
of the corresponding parameter:
\begin{itemize}
    \item \texttt{idB\&C-MILP}: solves IFD problem exactly as the
    MILP~\eqref{eqn:feasid} using an off-the-shelf solver;
    \item \texttt{idB\&C-MILP-k\_*}: solves IFD problem~\eqref{eqn:k-feasid}
    exactly for the specified value of $k$;
    \item \texttt{idB\&C-LS-k\_*}: employs the local search
    Algorithm~\ref{alg:enum} with the specified value of $k$; 
    \item \texttt{idB\&C-LS-k\_*-dBnd\_*\_*}: uses the local search
    Algorithm~\ref{alg:enum} when the depth of the  
    current node in the tree falls between the specified range. Otherwise, it
    solves~\eqref{eqn:feasid}; 
\end{itemize}
This last set of configurations tries to determine whether it is more
computationally advantageous to employ heuristics only at the lower levels of
the enumeration tree before resorting to solving~\eqref{eqn:feasid}, or vice
versa. Clearly, an exhaustive evaluation of all possible combinations of
different ranges and values of $k$ for each configuration would be
impractical. Therefore, based on the observations of Section~\ref{sec:id-gen},
we limited the range of \texttt{maxNeighborhoodSize} to be in $\brak{2, 3, 4,
5}$, and considered two settings for the ordered pair
\texttt{useLocalSearchDepthLb}\textbackslash\texttt{Ub}: $(0, 10)$ and $(10,
+\infty)$. For the case $k = 3$, we also tested additional pairs: $(0, 8)$,
$(0, 12)$, $(8, +\infty)$, and $(12, +\infty)$. In total, 21
different \texttt{idB\&C} configurations were tested.

As previously hinted, for those configurations using heuristics 
(i.e. all, except \texttt{idB\&C-MILP}) the solution 
of~\eqref{eqn:feasid} is mandatory if the current solution 
$(x, y) \in \S$ and no IFD is found by such approximations. 
This guarantees a correct application of Proposition~\ref{lemma:bfeas}. 

The configurations of \MIBS{} tested are as follows:
\begin{itemize}
    \item \MIBS{}: all parameters are set at their default values 
    (the generation of IDICs may be disabled); 
    \item \texttt{\MIBS{} only IDICs}: \emph{only} the generation of IDICs is
    enabled and no other classes; IFDs are generated by solving
    problem~\eqref{eqn:feasid} as an MILP using an off-the-shelf solver.  
    \item \texttt{\MIBS{} IDIC-MILP}: the generation of IDICs is enabled and
    separation of other classes of cuts is determined automatically by the
    default mechanism in \MIBS{}; IFDs are generated by solving
    problem~\eqref{eqn:feasid} as an MILP using an off-the-shelf solver.
    \item \texttt{\MIBS{} IDIC-LS-k\_*}: the generation of IDICs is enabled and
    separation of other classes of cuts is determined automatically by the
    default mechanism in \MIBS{}; IFDs are generated by first employing
    Algorithm~\ref{alg:enum} and then solving
    problem~\eqref{eqn:feasid} as an MILP when no IFD is found using local search.
\end{itemize}

The branch-and-cut scheme of \texttt{\MIBS{} only IDICs} differs 
from \texttt{idB\&C-MILP} solely in its use of the standard oracle. 
Therefore, it serves as the most relevant baseline configuration 
for comparison with \texttt{idB\&C}. In contrast, due to the use 
of multiple separation routines and the eventuality of disabling 
IDIC generation, \MIBS{} provides a fair comparison mainly with 
\texttt{\MIBS{} only IDICs}, \texttt{\MIBS{} IDIC-MILP}, and 
\texttt{\MIBS{} IDIC-LS-k\_*}. The specific values of $k$ and the control mechanism 
of the latter for finding IFDs were selected from the best-performing
configuration of \texttt{idB\&C}.

\subsection{Results}
To summarize the results of the experiment, we present several kinds of plots:
\begin{itemize}
    \item \textbf{Performance profiles} show empirical 
    cumulative distribution functions (CDFs) of ratios
    of a given performance measure of interest against the 
    \say{virtual best}~\citt{DolMorBenchmarking02}. Typical 
    measures are total solution time or the number of nodes 
    in the search tree. These plots are useful to compare the 
    performance of several configurations on the given measure. 

    \item \textbf{Baseline profiles} mimics the previous profiles, 
    but present empirical CDFs of ratios of a given performance 
    measure against the performance of a \say{baseline}, i.e. 
    a solver or a specific configuration of it, rather than 
    against the virtual best. These plots show more distinctly 
    the fraction of instances in which a certain configuration 
    outperforms (on the left side) or underperforms (on the right side) 
    the baseline.

    \item \textbf{Cumulative profiles} combine two profiles. 
    On the left side, they show the empirical CDFs of the 
    fraction of instances solved within the time limit. 
    On the right side, they report the empirical CDFs of the 
    fraction of instances that closed a certain final gap within 
    the time limit. Note that the lines on the two sides always 
    connect, since the fraction of instances solved within the 
    time limit equals that of the instances with zero gap at the time limit.
\end{itemize}

In order to create each plot on (a subset of) the dataset, we first solved all
instances with all considered configurations. Instances are then excluded from
the plots if (i) they were not solved within the time limit by any of the
configurations (except for cumulative profiles) and/or (ii) the
solution time is less than 5 seconds for all methods. 
Note that showing numerous configurations in a single plot may have a negative 
impact on its readability.
For this reason, each profile
will plot only the most promising results (highest performing
configurations). Those not shown can be considered less effective.

\begin{figure}[ht]
    \centering
    \begin{subfigure}{0.49\textwidth}
        \includegraphics[width=\linewidth]{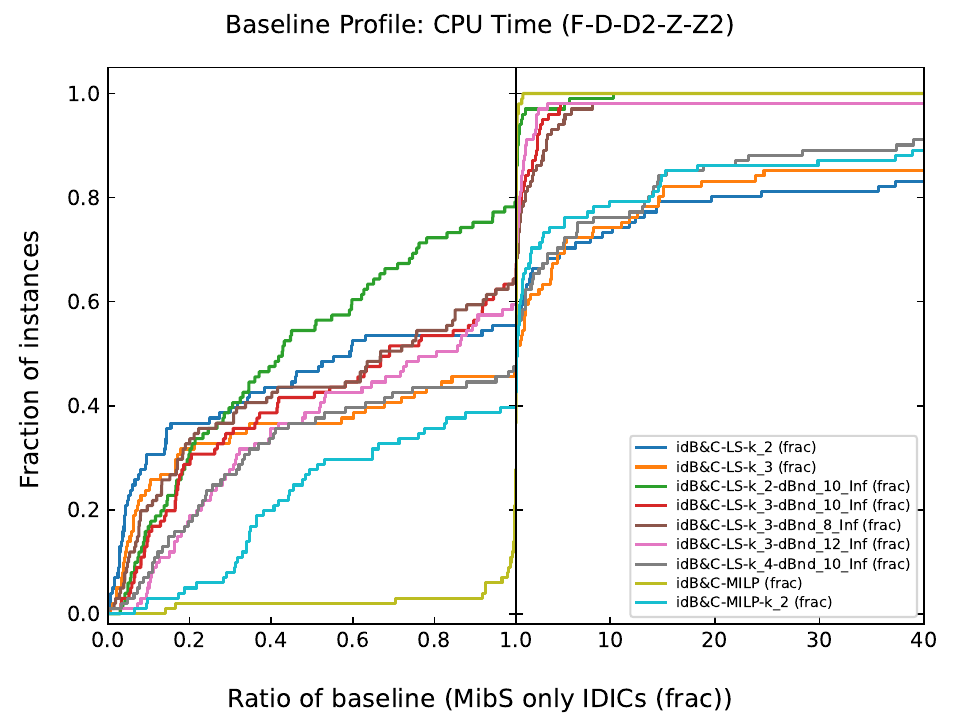}
        \caption{Baseline Profile Solution time}
        \label{fig:idbc:iblp:base:sol}
    \end{subfigure}    
    \begin{subfigure}{0.49\textwidth}
        \includegraphics[width=\linewidth]{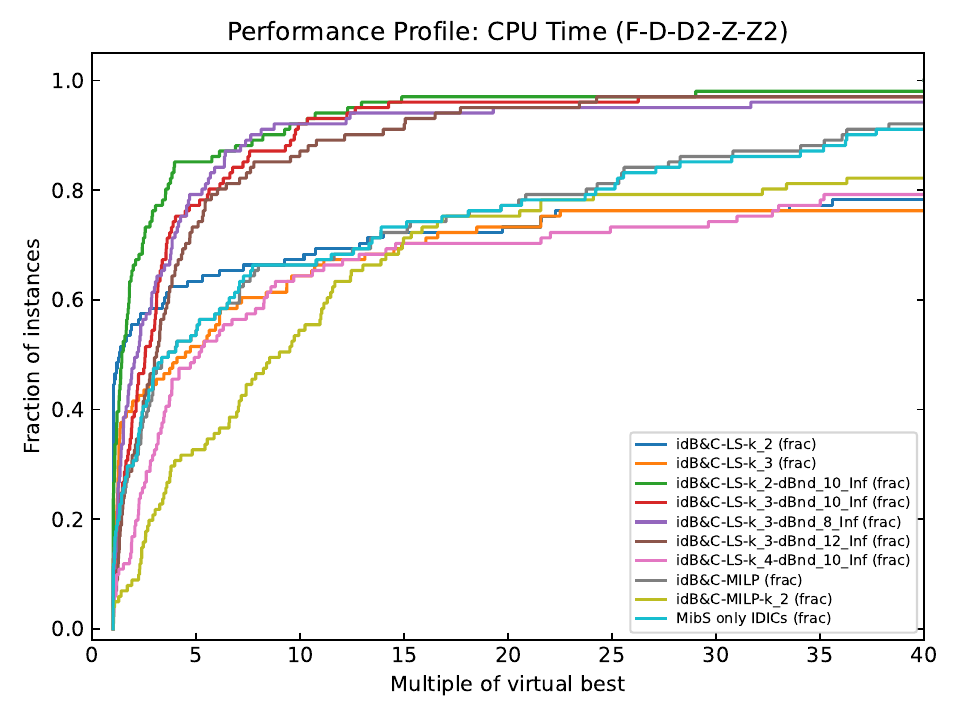}
        \caption{Performance Profile Solution time}
        \label{fig:idbc:iblp:perf:sol}
    \end{subfigure}
    \begin{subfigure}{0.49\textwidth}
        \includegraphics[width=\linewidth]{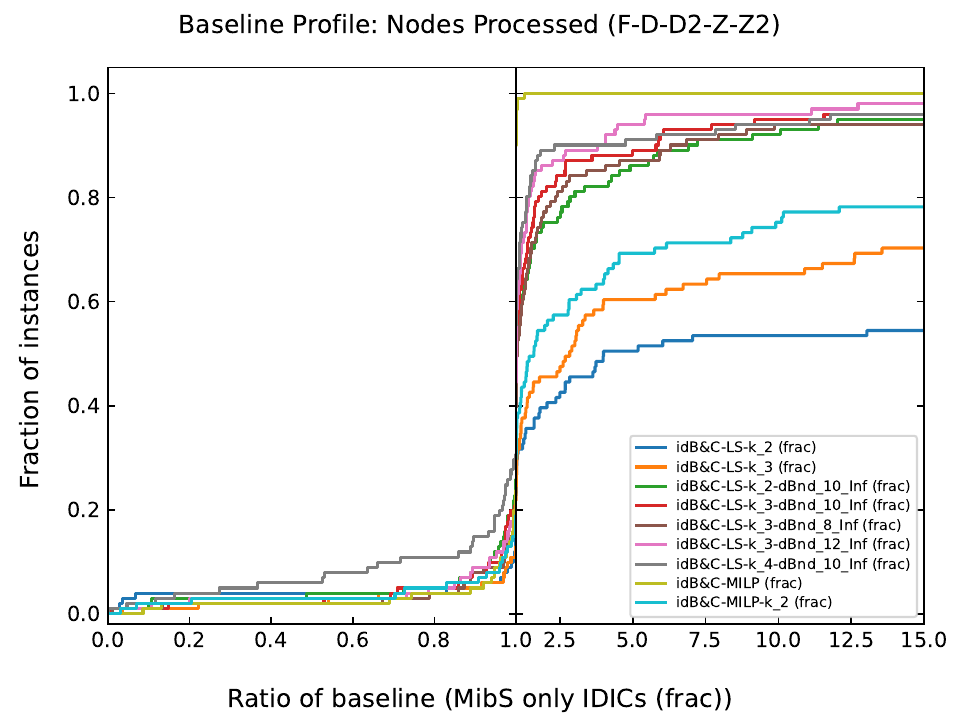}
        \caption{Baseline Profile nodes in the search tree}
        \label{fig:idbc:iblp:base:nodes}
    \end{subfigure}    
    \begin{subfigure}{0.49\textwidth}
        \includegraphics[width=\linewidth]{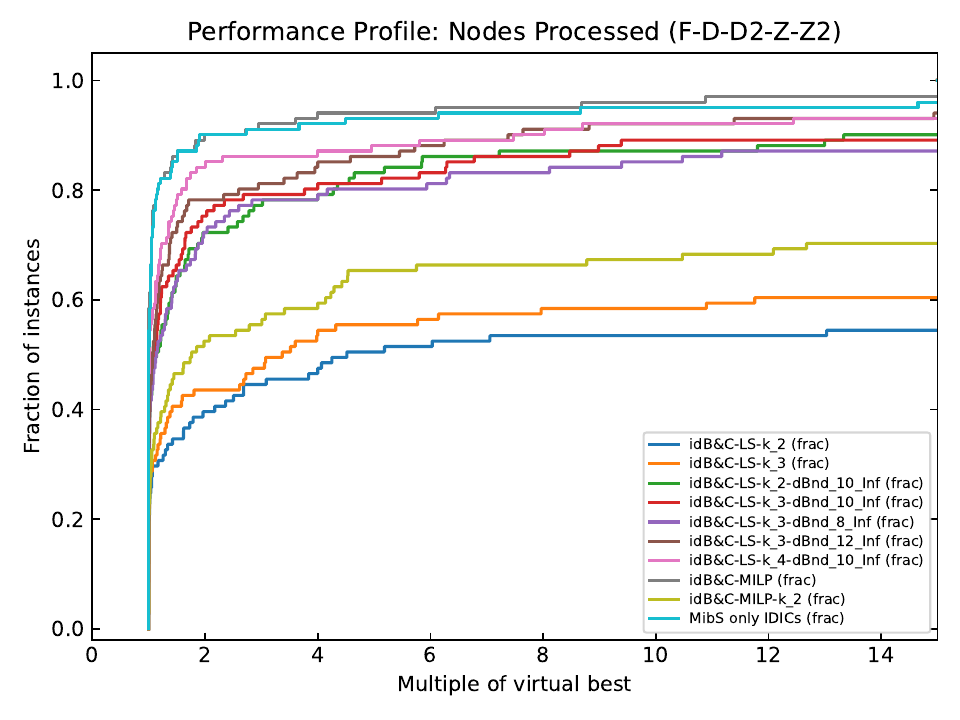}
        \caption{Performance Profile nodes in the search tree}
        \label{fig:idbc:iblp:perf:nodes}
    \end{subfigure} 
    \caption{Comparing \texttt{idB\&C} solution time and nodes in the search tree on IBLPs datasets}
    \label{fig:idbc:iblp:sol}   
\end{figure}

\begin{figure}
    \centering
    \begin{subfigure}{0.49\textwidth}
        \includegraphics[width=\linewidth]{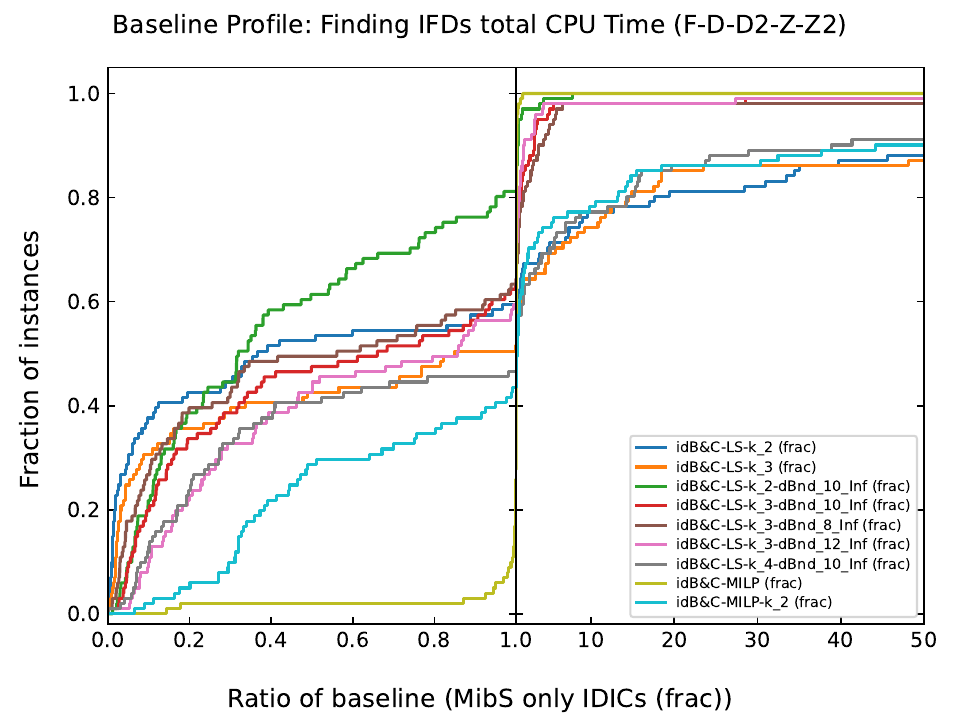}
        \caption{Baseline Profile Finding IFDs Tot. CPU time}
        \label{fig:idbc:iblp:base:ifd}
    \end{subfigure}    
    \begin{subfigure}{0.49\textwidth}
        \includegraphics[width=\linewidth]{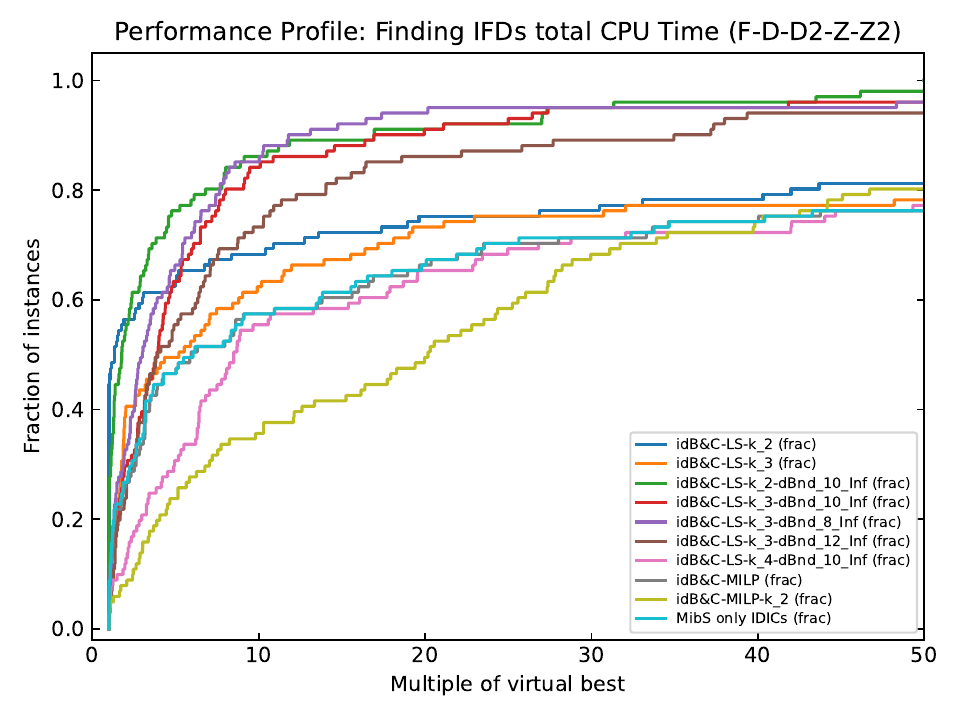}
        \caption{Performance Profile Finding IFDs Tot. CPU time}
        \label{fig:idbc:iblp:perf:ifd}
    \end{subfigure}  
    \begin{subfigure}{0.49\textwidth}
        \includegraphics[width=\linewidth]{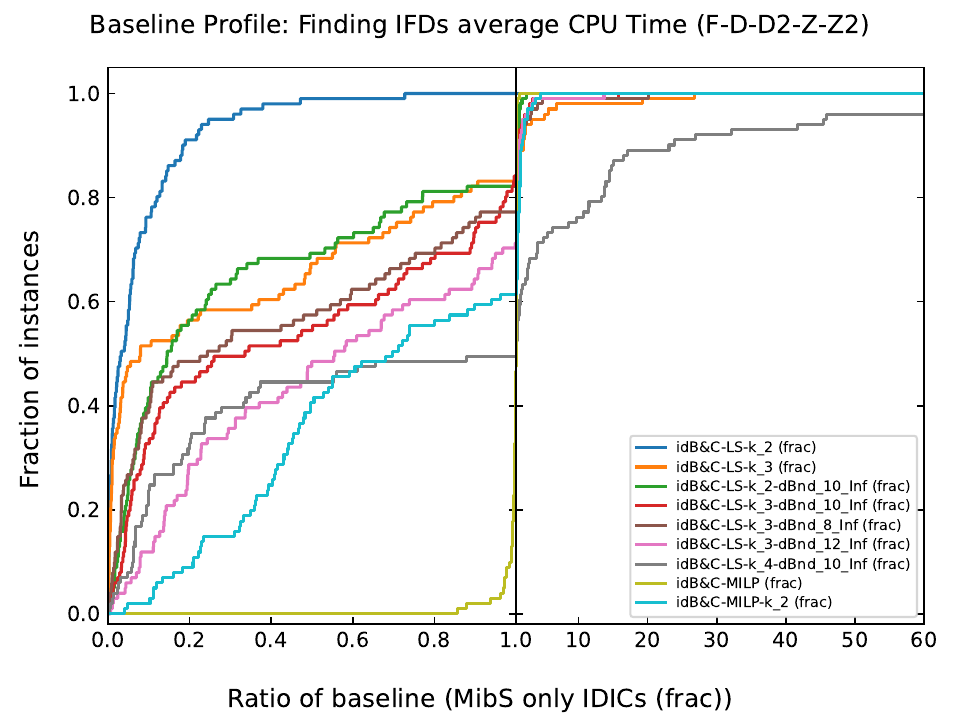}
        \caption{Baseline Profile Finding IFDs Avg. CPU time}
        \label{fig:idbc:iblp:base:avgifd}
    \end{subfigure}    
    \begin{subfigure}{0.49\textwidth}
        \includegraphics[width=\linewidth]{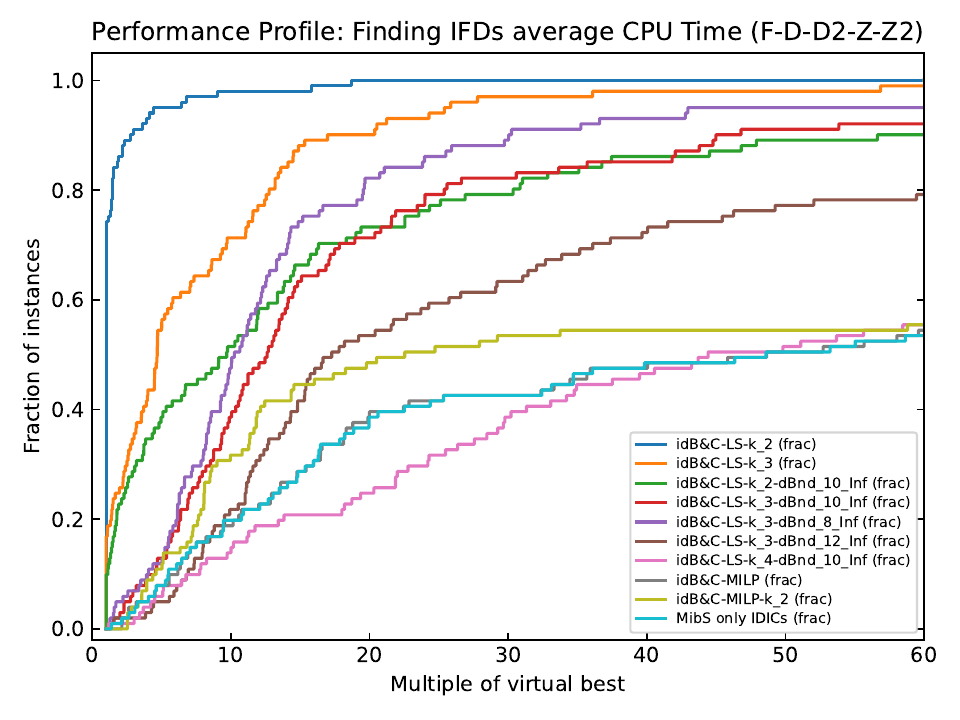}
        \caption{Performance Profile Finding IFDs Avg. CPU time}
        \label{fig:idbc:iblp:perf:avgifd}
    \end{subfigure}
    \caption{Comparing \texttt{idB\&C} times for finding IFDs on IBLPs datasets}
    \label{fig:idbc:iblp:ifd}       
\end{figure} 

\subsubsection{Improving Direction Oracle}

\begin{figure}
    \centering
    \begin{subfigure}{0.49\textwidth}
        \includegraphics[width=\linewidth]{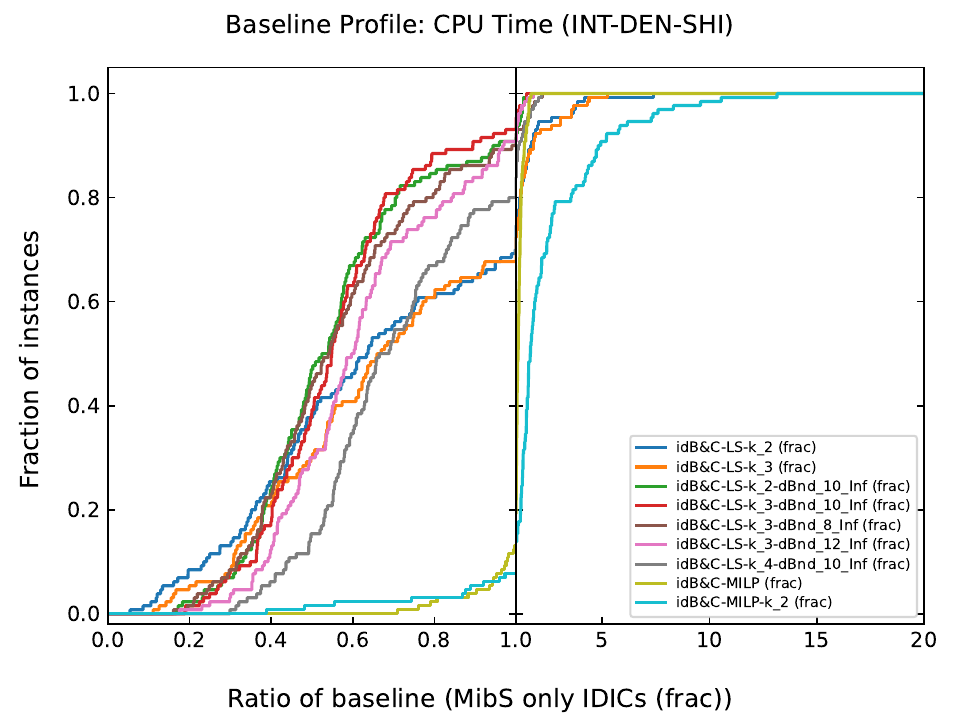}
        \caption{Baseline Profile Solution time}
        \label{fig:idbc:inter:base:sol}
    \end{subfigure}    
    \begin{subfigure}{0.49\textwidth}
        \includegraphics[width=\linewidth]{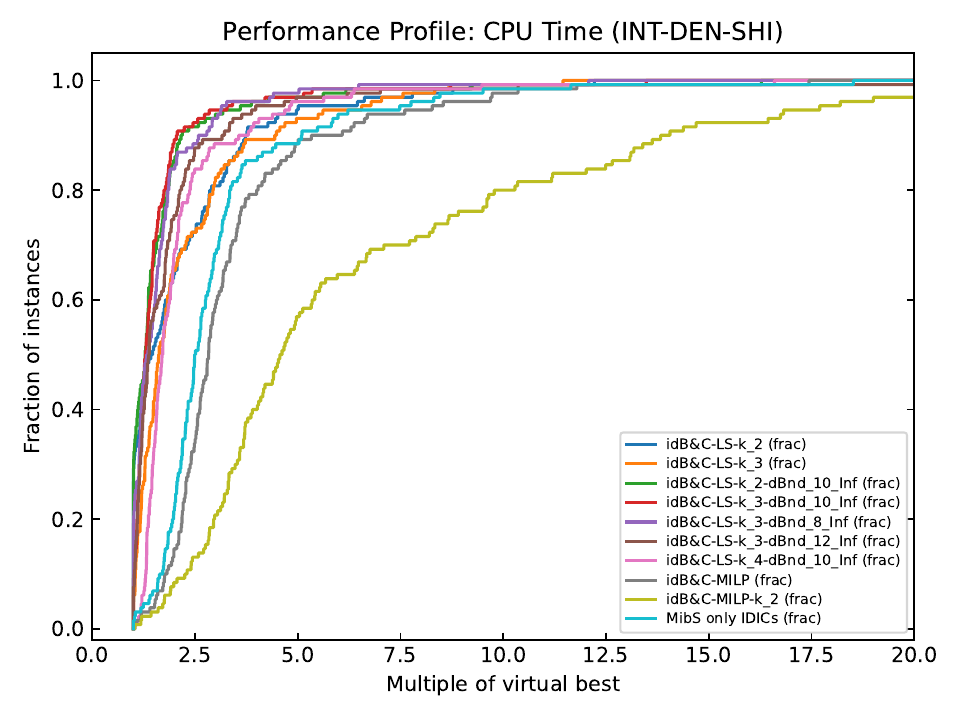}
        \caption{Performance Profile Solution time}
        \label{fig:idbc:inter:perf:sol}
    \end{subfigure} 
    \begin{subfigure}{0.49\textwidth}
        \includegraphics[width=\linewidth]{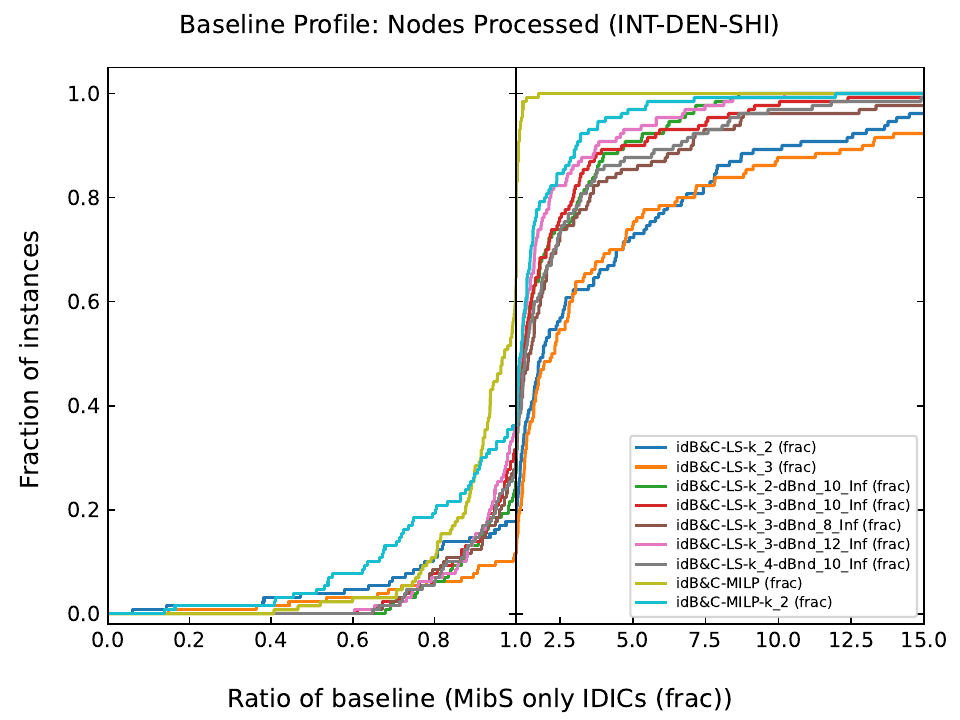}
        \caption{Baseline Profile nodes in the search tree}
        \label{fig:idbc:inter:base:nodes}
    \end{subfigure}    
    \begin{subfigure}{0.49\textwidth}
        \includegraphics[width=\linewidth]{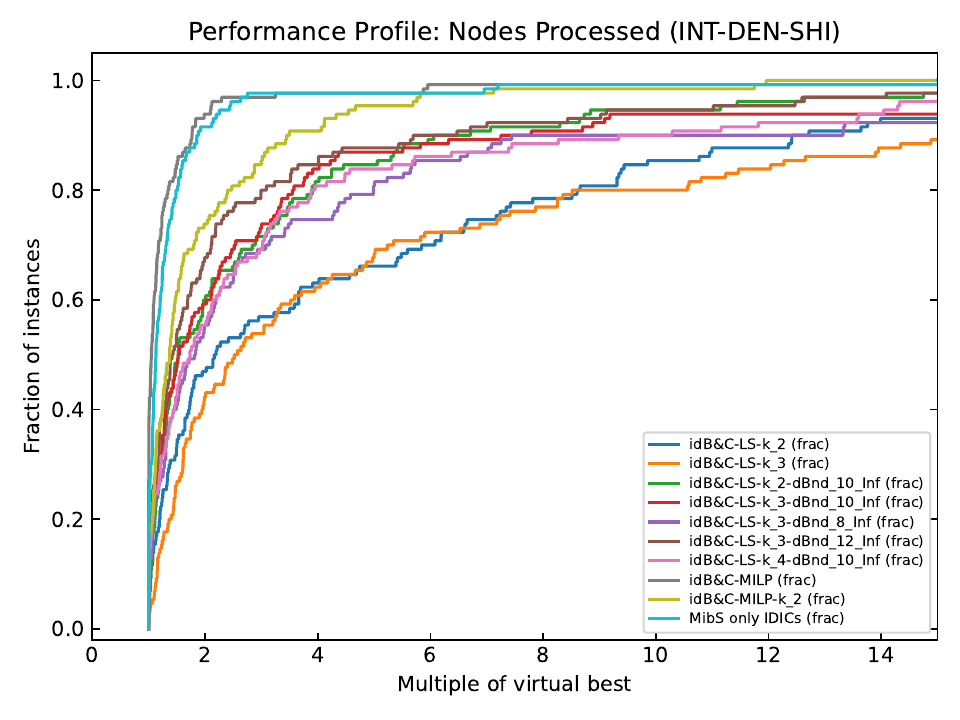}
        \caption{Performance Profile nodes in the search tree}
        \label{fig:idbc:inter:perf:nodes}
    \end{subfigure} 
    \caption{Comparing \texttt{idB\&C} solution time and nodes in the search tree on interdiction datasets}
    \label{fig:idbc:inter:sol}
\end{figure}

\begin{figure}
    \centering
    \begin{subfigure}{0.49\textwidth}
        \includegraphics[width=\linewidth]{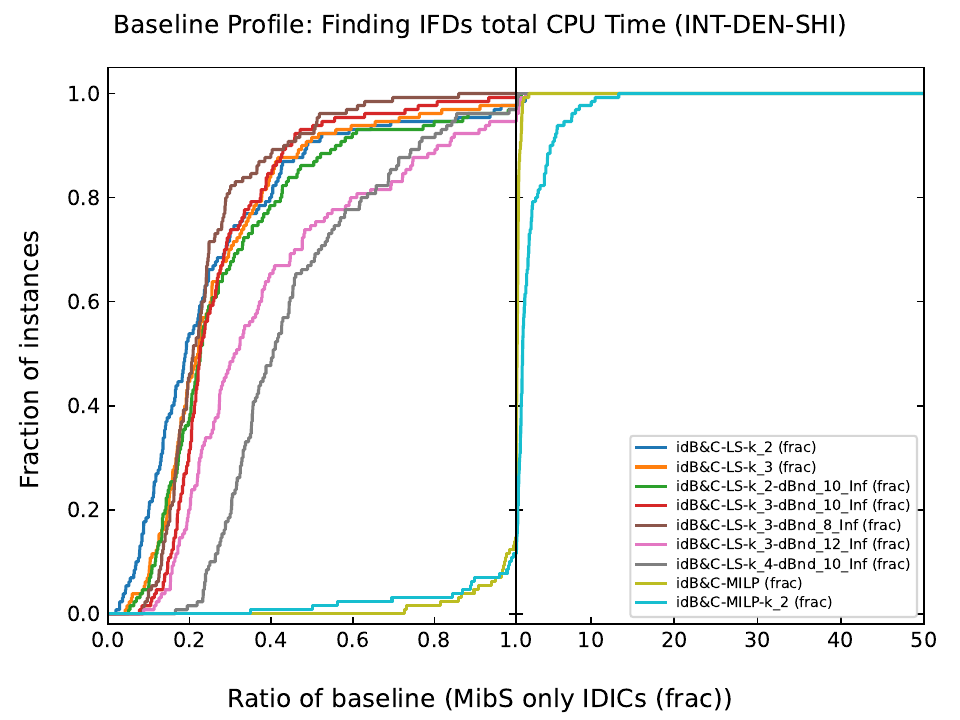}
        \caption{Baseline Profile Finding IFDs Tot. CPU time}
        \label{fig:idbc:inter:base:ifd}
    \end{subfigure}    
    \begin{subfigure}{0.49\textwidth}
        \includegraphics[width=\linewidth]{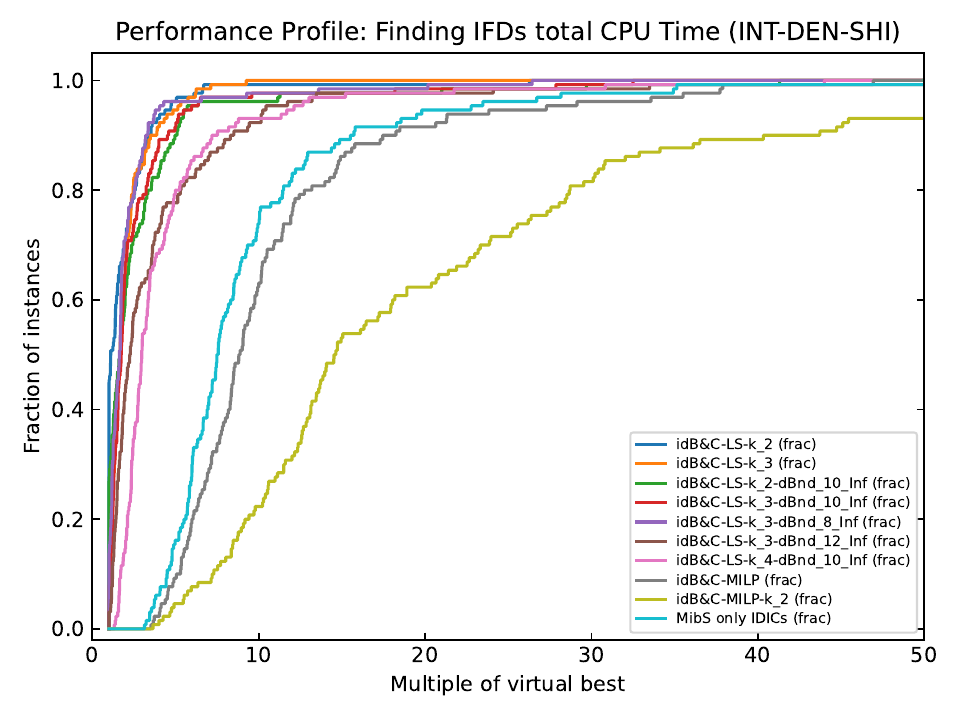}
        \caption{Performance Profile Finding IFDs Tot. CPU time}
        \label{fig:idbc:inter:perf:ifd}
    \end{subfigure}  
    \begin{subfigure}{0.49\textwidth}
        \includegraphics[width=\linewidth]{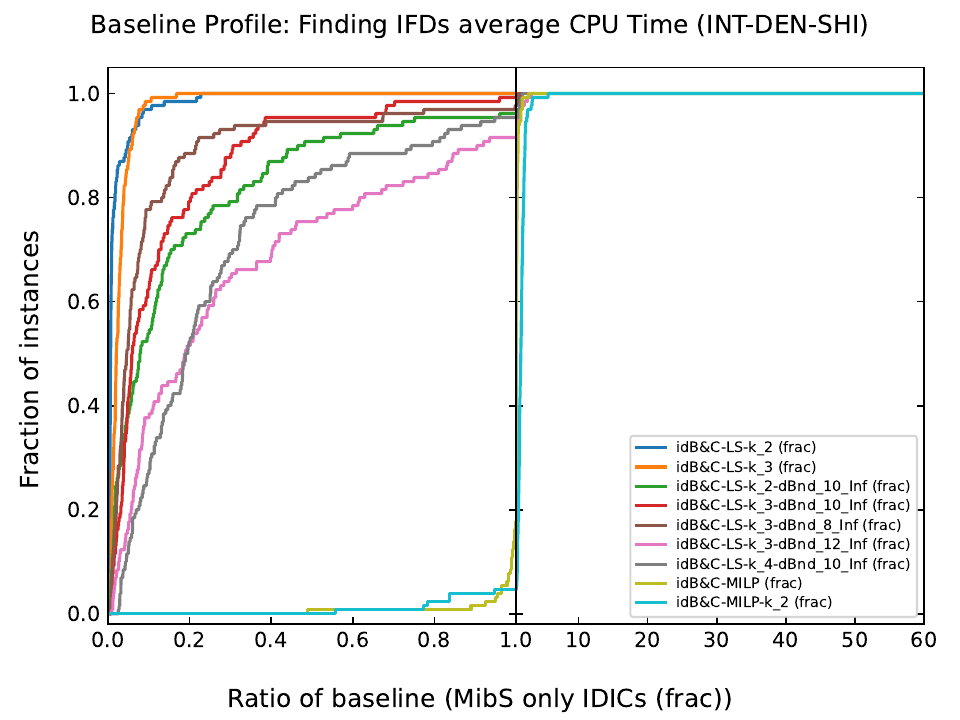}
        \caption{Baseline Profile Finding IFDs Avg. CPU time}
        \label{fig:idbc:inter:base:avgifd}
    \end{subfigure}    
    \begin{subfigure}{0.49\textwidth}
        \includegraphics[width=\linewidth]{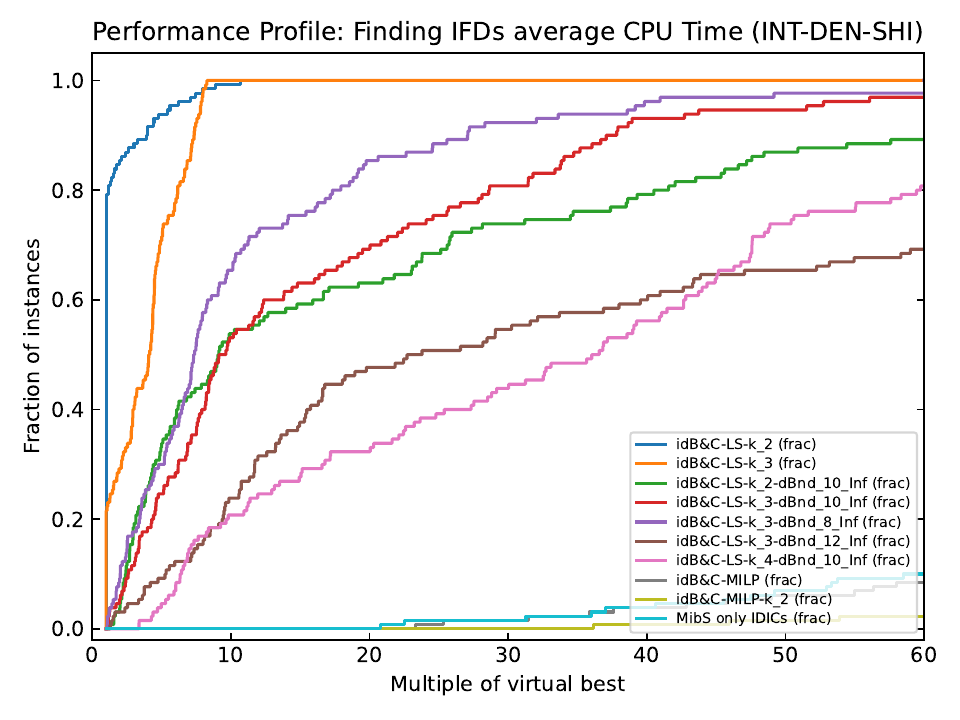}
        \caption{Performance Profile Finding IFDs Avg. CPU time}
        \label{fig:idbc:inter:perf:avgifd}
    \end{subfigure}    
    
    \caption{Comparing \texttt{idB\&C} times for finding IFDs on interdiction datasets}
    \label{fig:idbc:inter:ifd}  
\end{figure}

Figures~\ref{fig:idbc:iblp:sol}--\ref{fig:idbc:inter:ifd} show the results of
comparisons between \t{idB\&C} in various configurations with the \t{only
IDICs} configuration of \MIBS{}. This comparison provides the best
apples-to-apples comparison of \MIBS{} and the algorithmic framework proposed
here, holding as many things constant as possible. We display baseline and
performance profiles using solution CPU time, number of nodes in the search
tree, total and average CPU time spent in finding IFDs as measures on all IBLP
and Interdiction datasets, comparing all \texttt{idB\&C} configurations
against \texttt{\MIBS{} only IDICs} (used as the baseline).

The first remarkable result deduced from this experiment is that
\texttt{idB\&C-MILP} shows very similar performances to \texttt{\MIBS{} only IDICs}
over both classes of instances. Recall that \texttt{idB\&C} never solves the
follower's problem explicitly, but instead exploits only
Proposition~\ref{lemma:bfeas} to detect bilevel feasible solutions. Thus, this
change by itself does not yield any improvement, but it also does not degrade
performance.

From Figure~\ref{fig:idbc:iblp:base:nodes}
and~\ref{fig:idbc:inter:base:nodes}, we observe the expected tradeoff between
the time spent generating a ``good'' IFD versus the strength of the resulting
inequalities. The results show that employing local
search for finding IFDs generally comes at a cost in terms of the strength of
the separated IDICs. This is expected when replacing an exact method with a
heuristic and the result is an increase in the search tree size.
In terms of total CPU time, as well as time spent finding IFDs, however, the
trade-off strongly favors the use of local search.
As a matter of fact,
Figures~\ref{fig:idbc:iblp:base:sol} and~\ref{fig:idbc:iblp:base:ifd}
(\ref{fig:idbc:inter:base:sol} and~\ref{fig:idbc:inter:base:ifd}, resp.)
indicate that local search results in a significant reduction in both the time
spent to identify IFDs and total CPU time on a substantial fraction of IBLP
(Interdiction, resp.) instances, across all configuration. 

In this regard, the strategy of solving~\eqref{eqn:feasid} exactly only in
higher levels of the tree is key in balancing the drawbacks due to the
employment of a heuristic approach with advantages. Interestingly, the
experiment emphasize that the best trade-off is achieved for small values of
$k$, i.e., 2 and 3. This result corroborates the related work
of~\citt{XueProRal22}.

More specifically, consider the solution CPU time on the IBLP dataset.
Figures~\ref{fig:idbc:iblp:base:sol} and~\ref{fig:idbc:iblp:perf:sol} show
that the configurations \texttt{idB\&C-LS-k\_*-dBnd\_*\_Inf} are the
best-performing settings. While both \texttt{idB\&C-LS-k\_*}
and \texttt{idB\&C-MILP-k\_2} reduce the solution time on a non-trivial
portion of the dataset compared to the baseline, they also exhibit significant
slowdowns on the remaining fraction. In
contrast, \texttt{idB\&C-LS-k\_2-dBnd\_10\_Inf} offers more consistent
performance, achieving at least a 20\% improvement over the baseline for
roughly 75\% of instances, which increases to at least 50\% for over half the
dataset.

An even more prominent result is highlighted in
Figure~\ref{fig:idbc:inter:base:sol} on the Interdiction dataset, where
Algorithm~\ref{alg:enum} is particularly effective in finding IFDs on problem
with such structure. The more pronounced gains are again obtained
by \texttt{idB\&C-LS-k\_*-dBnd\_10\_Inf}, for $k \in \brak{2, 3}$, showing a
significant reduction in solution time over more than 90\% of all the dataset,
and similar time to the baseline on the remaining fraction.

\subsubsection{Using Local Search in \MIBS{}}

In the next experiment, we compares all previously described configurations
of \MIBS{} in order to measure the effect of using the local search for
finding IFDs, holding other elements of the algorithmic strategy of \MIBS{}
constant. Given the promising outcome of the previous experiment, the
parameter settings and the value of $k$ for the best-performing configuration
of \t{idB\&C} (\t{idB\&C-LS-k\_2-dBnd\_10\_Inf}) were replicated for the
configuration of \MIBS{} using local search, i.e., \texttt{\MIBS{} IDIC-LS-k\_2}.

\begin{figure}
    \centering
    \begin{subfigure}{0.49\textwidth}
        \includegraphics[width=\linewidth]{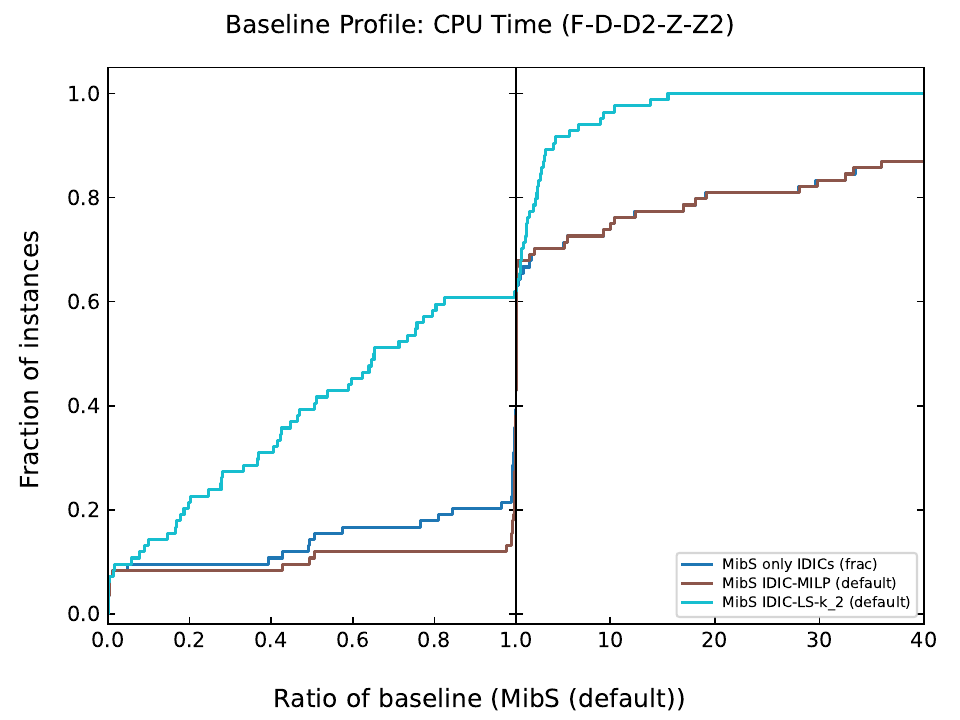}
        \caption{Baseline Profile Solution time}
        \label{fig:mibs:iblp:base:sol}
    \end{subfigure}    
    \begin{subfigure}{0.49\textwidth}
        \includegraphics[width=\linewidth]{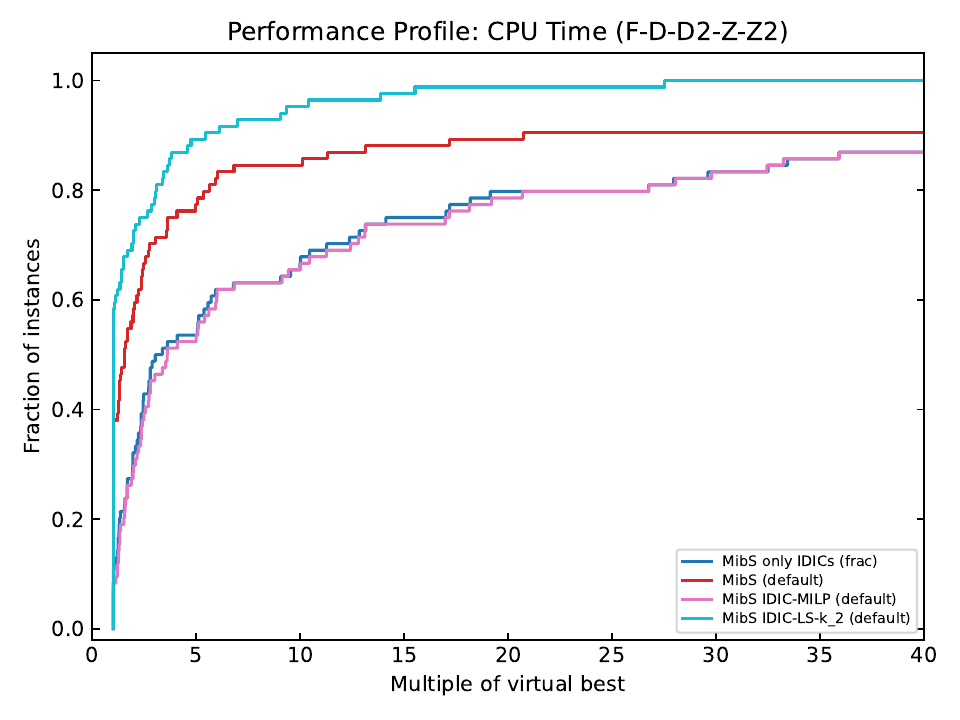}
        \caption{Performance Profile Solution time}
        \label{fig:mibs:iblp:perf:sol}
    \end{subfigure} 
    \begin{subfigure}{0.49\textwidth}
        \includegraphics[width=\linewidth]{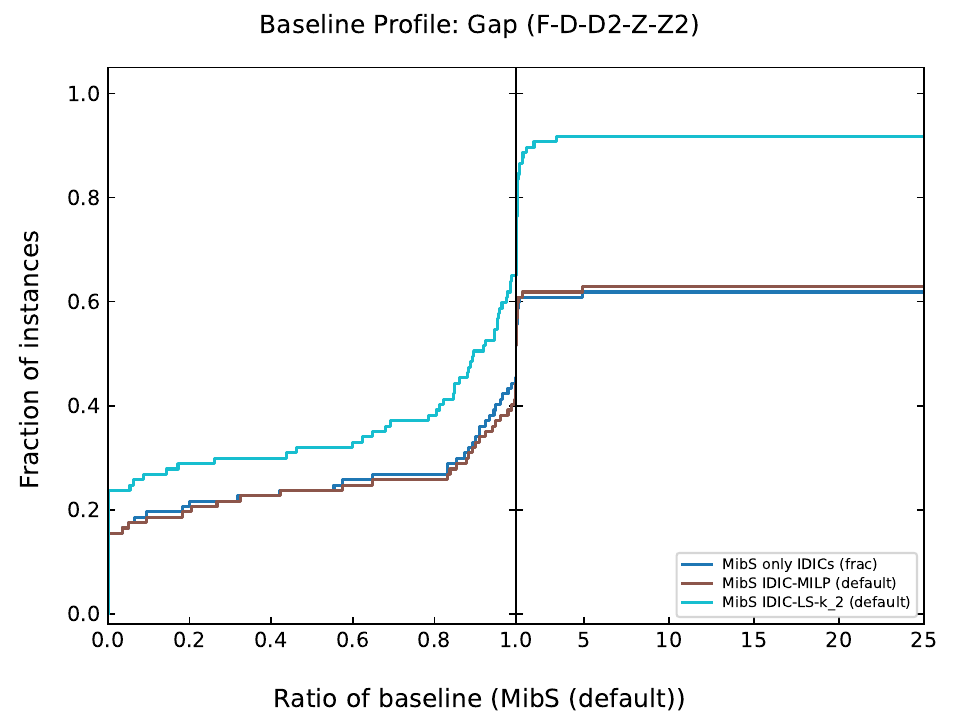}
        \caption{Baseline Profile Gap}
        \label{fig:mibs:iblp:base:gap}
    \end{subfigure}      
    \begin{subfigure}{0.49\textwidth}
        \includegraphics[width=\linewidth]{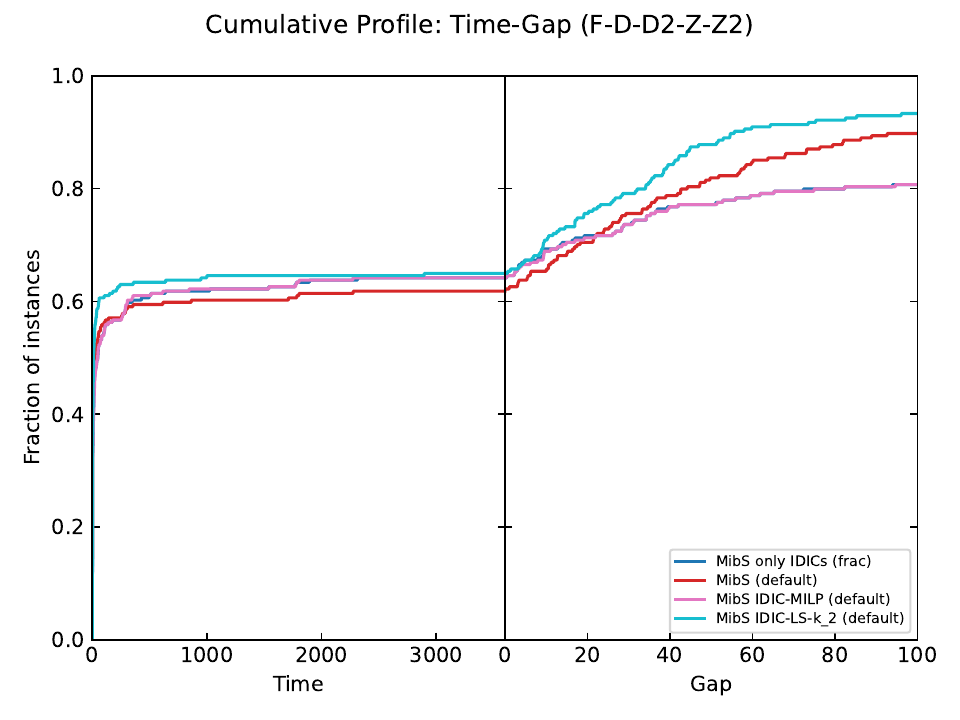}
        \caption{Cumulative Profile}
        \label{fig:mibs:iblp:cum}
    \end{subfigure} 
    \begin{subfigure}{0.49\textwidth}
        \includegraphics[width=\linewidth]{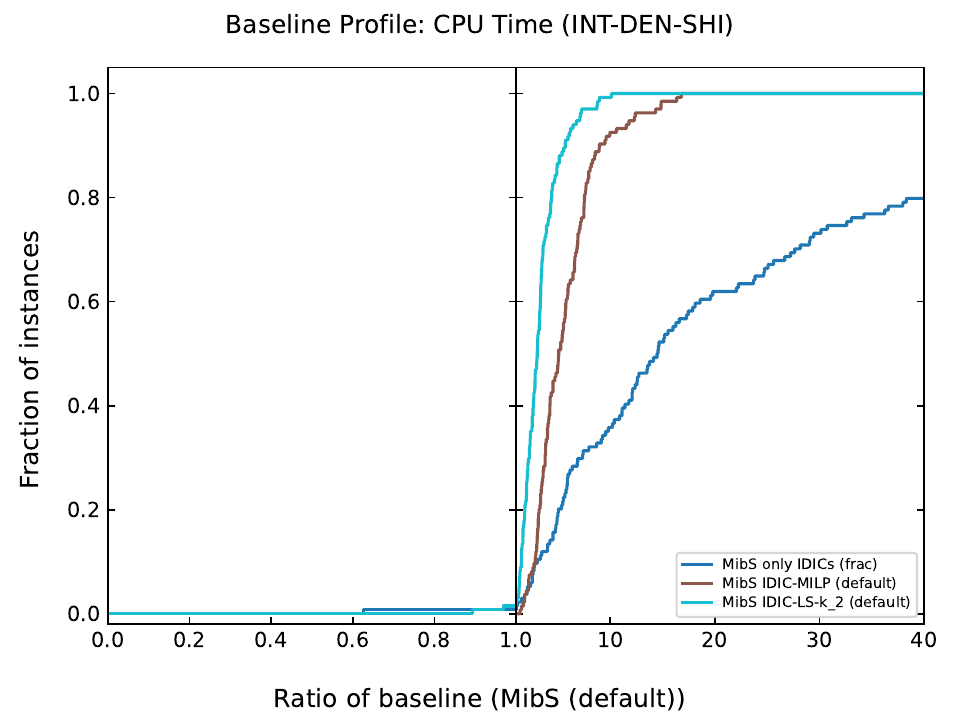}
        \caption{Baseline Profile Solution time}
        \label{fig:mibs:inter:base:sol}
    \end{subfigure}    
    \begin{subfigure}{0.49\textwidth}
        \includegraphics[width=\linewidth]{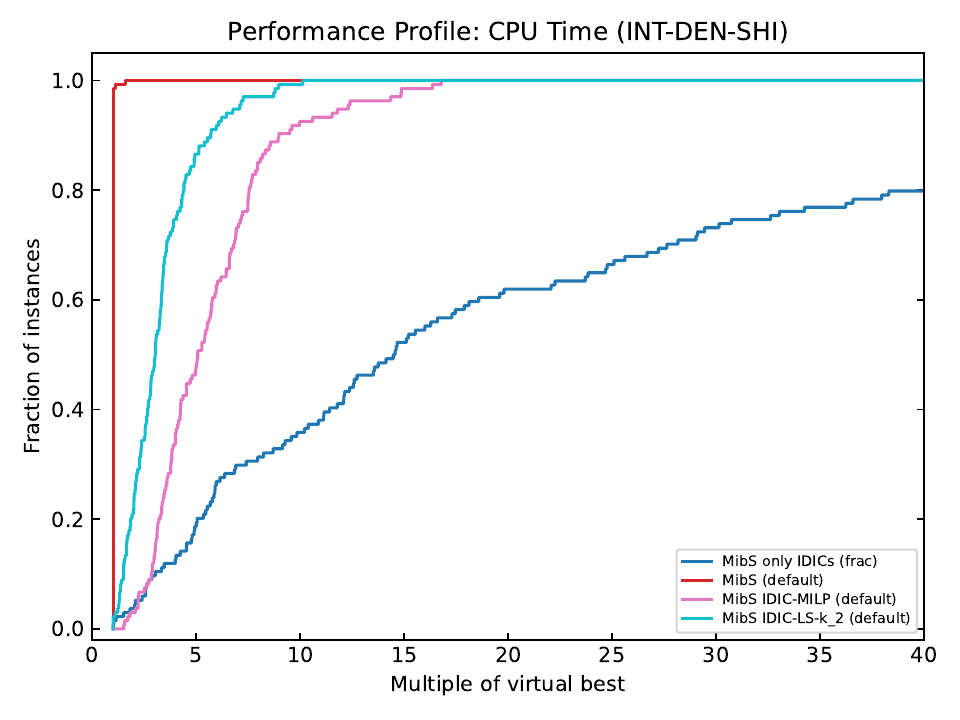}
        \caption{Performance Profile Solution time}
        \label{fig:mibs:inter:perf:sol}
    \end{subfigure}    
    
    \caption{Comparing \MIBS{} solution times on IBLP and Interdiction datasets}
    \label{fig:mibs:all}
\end{figure}

Figure~\ref{fig:mibs:all} illustrates baseline, performance and 
cumulative profiles using CPU time and final gap closed as performance
measures, with the default \MIBS{} configuration serving as the baseline.  
The CPU time spent in finding IFDs is excluded since \MIBS{} might generate
different classes of cuts and IDICs may be disabled.

On the one hand, we note from Figure~\ref{fig:mibs:inter:base:sol} that \MIBS{} is the
clear winner on the interdiction dataset, primarily due to the use of cuts
specialized for these problems. Nevertheless,
Figure~\ref{fig:mibs:inter:perf:sol} shows that \texttt{\MIBS{} IDIC-LS-k\_2}
is the best performing configuration among all other instances.
Figure~\ref{fig:mibs:iblp:base:sol} shows that on IBLPs, \texttt{\MIBS{} IDIC-LS-k\_2}
outperforms \MIBS{} by at least 20\% on 60\% of instances, with solution time
reductions reaching up to 60\% on about 30\% of the dataset. Importantly, by
examining the performance of \texttt{\MIBS{} IDIC-MILP} one can see that this
improvement is attributed not to the use of IDICs alone but to the application
of local search specifically. Furthermore,
Figures~\ref{fig:mibs:iblp:base:gap} and~\ref{fig:mibs:iblp:cum} indicate
that \MIBS{} is able to close more gap when local search is used. In
particular, \texttt{\MIBS{} IDIC-LS-k\_2} is able to solve about 12 additional
instances to optimality compared to \texttt{\MIBS{} IDIC-MILP} and more than
50 with respect to \MIBS{}.

\begin{figure}
    \centering
    \begin{subfigure}{0.49\textwidth}
        \includegraphics[width=\linewidth]{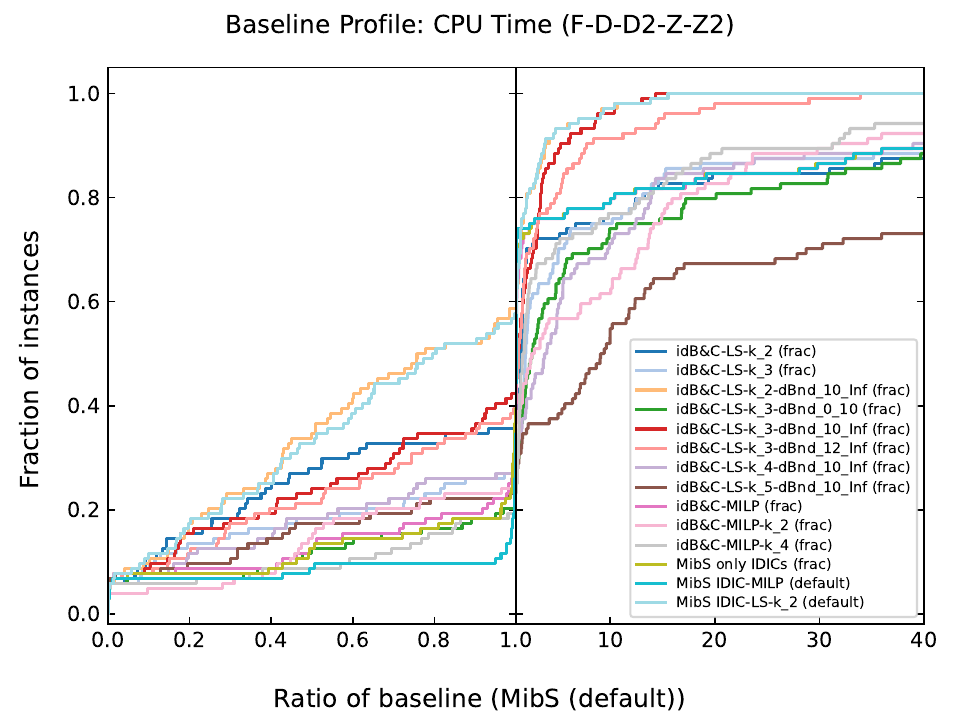}
        \caption{Baseline Profile Solution time}
        \label{fig:all:iblp:base:sol}
    \end{subfigure}    
    \begin{subfigure}{0.49\textwidth}
        \includegraphics[width=\linewidth]{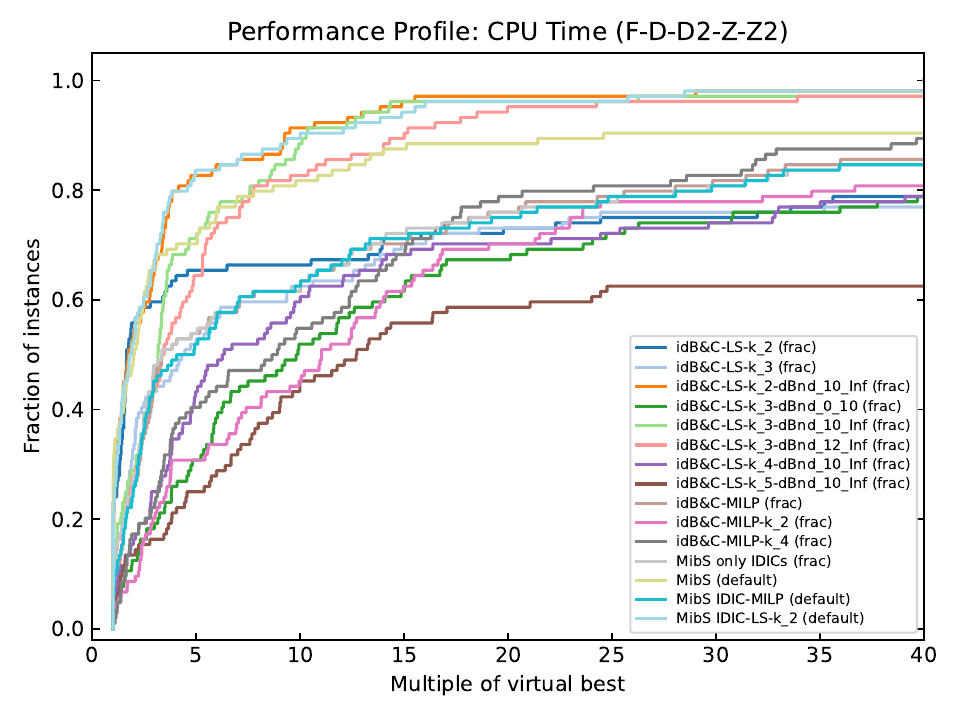}
        \caption{Performance Profile Solution time}
        \label{fig:all:iblp:perf:sol}
    \end{subfigure}     
    \begin{subfigure}{0.49\textwidth}
        \includegraphics[width=\linewidth]{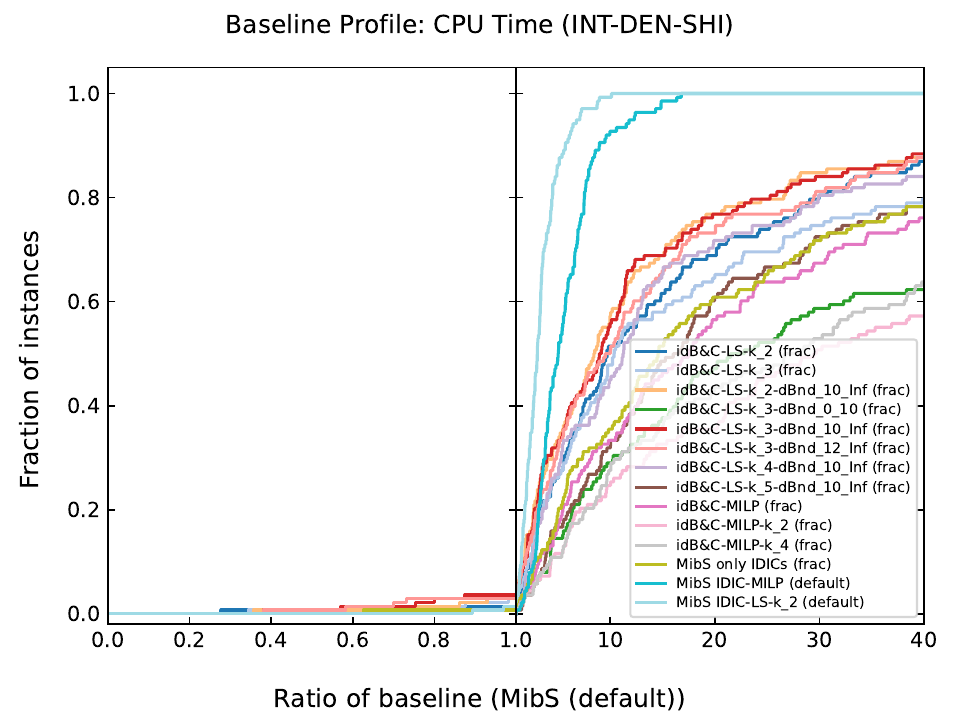}
        \caption{Baseline Profile Solution time}
        \label{fig:all:inter:base:sol}
    \end{subfigure}    
    \begin{subfigure}{0.49\textwidth}
        \includegraphics[width=\linewidth]{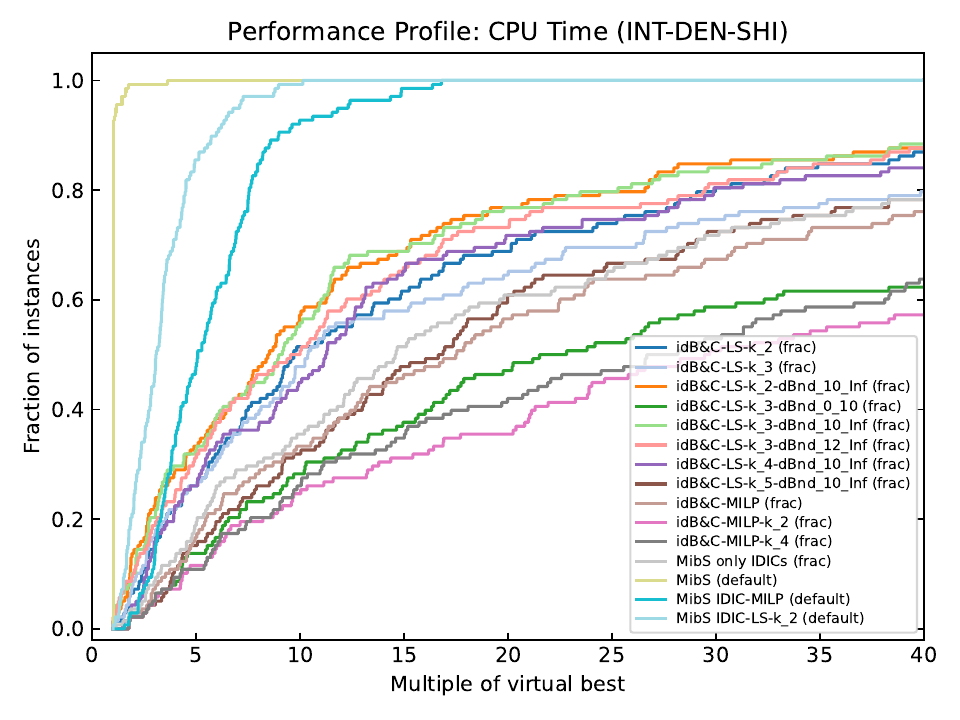}
        \caption{Performance Profile Solution time}
        \label{fig:all:inter:perf:sol}
    \end{subfigure}      
    
    \caption{Comparing all configurations on all the datasets}
    \label{fig:all}
\end{figure}

Finally, Figure~\ref{fig:all} shows baseline and performance profiles plotting
the majority of configurations tested in this analysis. Interestingly, 
plots~\ref{fig:all:iblp:base:sol} and~\ref{fig:all:iblp:perf:sol} reveal that 
\texttt{\MIBS{} IDIC-LS-k\_2} and \texttt{idB\&C-LS-k\_2-dBnd\_10\_Inf} have
very similar performance, emphasizing the competitiveness of the oracle based
on improving directions.

\section{Conclusions}
\label{sec:conclusions}

Improving directions are a fundamental and versatile tool in mixed-integer bilevel 
linear optimization, underpinning branching schemes, the formulation of 
optimality-based relaxations and the generation of strong valid inequalities, 
and they have been instrumental in advancing the state-of-the-art solution 
methods for this class of problems. 

In this work, we have taken an important step toward explaining and quantifying how 
improving directions contribute to restoring the follower's optimality condition, 
by showing that the convex hulls of the feasible regions 
arising from the optimality-based hierarchy of relaxations are exactly characterized
from valid inequalities stemming from improving directions.
Moreover, we have shed light on a new role that improving directions play, as they 
unify the oracle computations for checking bilevel feasibility and generating 
strong valid inequalities, a perspective that may lead to substantial improvements 
in empirical performance, as suggested by our promising experiments with our
branch-and-cut framework. 

Notably, this novel algorithm is not limited to the generation of intersection cuts 
discussed here, but can be easily extended to any class of valid inequality for MIBLPs
somehow related to the existence of an improving solution or direction.
Since our branch-and-cut highlights the centrality of finding an improving direction
as \NPcomplexity-complete subproblem,
we plan to explore the integration of MILP warm-start capabilities implemented 
in SYMPHONY~\citp{SYMPHONY5.7} (and already used in \MIBS{} as subsolver), 
to further enhance performance.

Due to the promising results of the local search
shown here, we also plan to develop more refined mechanism for 
better controlling its dynamic employment and other related enhancements.   
From a theoretical perspective instead, we will explore 
connection with more general, yet related, optimality certificates given by 
test sets for pure integer linear problems as already investigated in, 
e.g.,~\citp{ConTraBuchberger91,ScaTest97}.

{\footnotesize \textbf{Acknowledgements}\quad This research was made possible with support from Office of Naval Research Grant N000142212676.}

\clearpage
\pagebreak


\end{document}